\documentclass[11pt,leqno]{amsart}
\usepackage[frenchb]{babel}
\usepackage{amsmath,amsthm,amsfonts,amssymb,amscd}

\newtheorem{theo}{Th\'eor\`eme}
\newtheorem{lem}{Lemme}
\newtheorem{defi}{D\'efinition}

\newcommand{\al}{\alpha}
\newcommand{\bt}{\beta}
\newcommand{\ga}{\gamma}
\newcommand{\dt}{\delta}
\newcommand{\Dt}{\Delta}
\newcommand{\ep}{\epsilon}
\newcommand{\la}{\lambda}
\newcommand{\La}{\Lambda}
\newcommand{\ph}{\varphi}
\newcommand{\ta}{\vartheta}
\newcommand{\si}{\sigma}
\newcommand{\Si}{\Sigma}
\newcommand{\om}{\omega}
\newcommand{\G}{\mathcal{G}}
\newcommand{\K}{\mathcal{K}}
\newcommand{\R}{\mathcal{R}}
\newcommand{\T}{\mathcal{T}}
\newcommand{\U}{\mathcal{U}}
\newcommand{\V}{\mathcal{V}}
\newcommand{\W}{\mathcal{W}}
\newcommand{\BS}{\mathbf{S}}
\newcommand{\ti}{\times}
\newcommand{\su}{\subset}
\newcommand{\sm}{\setminus}
\newcommand{\ii}{\infty}
\newcommand{\vi}{\emptyset}
\newcommand{\ov}{\overline}
\newcommand{\td}{\tilde}
\newcommand{\wt}{\widetilde}
\newcommand{\wh}{\widehat}
\newcommand{\mt}{\multimap}
\newcommand{\pa}{\partial}
\newcommand{\oc}[1]{\overset{\circ}#1}

\DeclareMathOperator{\Tr}{Tr} \DeclareMathOperator{\St}{St}
\DeclareMathOperator{\diam}{diam} \DeclareMathOperator{\Fix}{Fix}
\DeclareMathOperator{\ind}{ind} \DeclareMathOperator{\im}{im}

\begin{document}
\title[Points fixes des applications compactes]{Points fixes des applications compactes
dans les espaces ULC}
\author{Robert Cauty}
\address{Universit\'e Paris 6, Institut de Math\'ematiques de Jussieu, case 247, 4 place
Jussieu, 75252 Paris Cedex 05} \email{cauty@math.jussieu.fr}

\begin{abstract} A topological space is locally equiconnected if there exists a
neighborhood $U$ of the diagonal in $X\ti X$ and a continuous map $\la:U\ti [0,1]\to X$
such that $\la(x,y,0)=x$, $\la(x,y,1)=y$ and $\la(x,x,t)=x$ for $(x,y)\in U$ and
$(x,t)\in X\ti[0,1]$. This class contains all ANRs, all locally contractible topological
groups and the open subsets of convex subsets of linear topological spaces. In a series
of papers, we extended the fixed point theory of compact continuous maps, which was well
developped for ANRs, to all separeted locally equiconneced spaces. This generalization
includes a proof of Schauder's conjecture for compact maps of convex sets. This paper is
a survey of that work. The generalization has two steps, which needs entirely different
methods: the metrizable case, and the passage from the metrizable case to the general
case. The metrizable case is, by far, the most difficult. To treat this case, we
introduced in [4] the notion of algebraic ANR. Since the proof that metrizable locally
equiconnected spaces are algebraic ANRs is rather difficult, we give here a detailed
sketch of it in the case of a compact convex subset of a metrizable t.v.s.. The passage
from the metrizable to the general case uses a free functor and representations of
compact spaces as inverse limits of some special inverse systems of metrizable compacta.
\end{abstract}

\maketitle
\section{Introduction}

Soit $I=[0,1]$. Un espace topologique $X$ est dit uniform\'ement localement contractile,
ou ULC, s'il existe un voisinage $U$ de la diagonale de $X\ti X$ et une fonction continue
$\la:U\ti I\to X$ v\'erifiant $\la(x,y,0)=x$, $\la(x,y,1)=y$ et $\la(x,x,t)=x$ pour
$(x,y)\in U$ et $(x,t)\in X\ti I$. Si une telle fonction existe sur $X\ti X\ti I$, alors
$X$ est dit uniform\'ement contractile, ou UC. La classe des espaces ULC contient tous
les sous-ensembles convexes des e.v.t., tous les groupes topologiques localement
contractiles et tous les r\'etractes de voisinages de tels espaces. En particulier, tout
r\'etracte absolu de voisinage est ULC, mais l'exemple construit dans [1] montre que les
espaces ULC m\'etrisables ne sont pas tous des r\'etractes absolus de voisinage.

La th\'eorie des points fixes des applications compactes (univoques ou multivoques) est
bien d\'evelopp\'ee dans la classe des r\'etractes absolus de voisinage et dans celle des
e.v.t. localement convexes. Malheureusement, les d\'emonstrations utilis\'ees dans ces
deux classes ne s'appliquent pas \`a tous les espaces m\'etriques lin\'eaires, car elles
emploient des propri\'et\'es de prolongement ou d'approximation des applications
compactes que, comme le montre l'exemple de [1], ne poss\`edent pas tous les espaces
m\'etriques lin\'eaires. Il est cependant possible d'\'etendre aux espaces ULC tous les
th\'eor\`emes de point fixe connus pour les applications compactes dans les r\'etractes
absolus de voisinage, et le but de cet article est d'expliquer comment faire
syst\'ematiquement cette extension.

Selon un sch\'ema introduit dans [2], la d\'emonstration d'un r\'esultat de point fixe
pour les applications compactes dans les espaces ULC se fait en deux \'etapes, qui
n\'ecessitent des techniques compl\`etement diff\'erentes: le cas des espaces
m\'etrisables et le passage du cas m\'etrisable au cas g\'en\'eral.

Le cas m\'etrisable est la partie la plus d\'elicate. L'approche utilis\'ee dans [2] est
\`a oublier, \'etant avantageusement supplant\'ee par la th\'eorie des r\'etractes
absolus de voisinage alg\'ebriques d\'evelopp\'ee dans [4], qui fournit des r\'esultats
beaucoup plus forts (il y a d'ailleurs une erreur dans la d\'emonstration du lemme 3 de
[2], qu'il n'y a plus de raison de corriger, vu la sup\'eriorit\'e de la nouvelle
approche). La th\'eorie des r\'etractes absolus de voisinage alg\'ebriques permet aussi
de comprendre pourquoi les espaces m\'etriques lin\'eaires ont la propri\'et\'e du point
fixe pour les applications compactes alors que ce ne sont pas n\'ecessairement des
extenseurs pour la classe des compacts m\'etrisables : c'est une manifestation des
diff\'erences de propri\'et\'es entre les groupes d'homologie et ceux d'homotopie ! La
d\'efinition et les principales propri\'et\'es des r\'etractes absolus de voisinage
alg\'ebriques seront rappel\'ees \`a la section 2. Dans les r\'etractes absolus de
voisinage alg\'ebriques, il est naturel d'\'etudier les points fixes des fonctions
continues par l'interm\'ediaire des morphismes qu'elles induisent sur le complexe des
cha\^\i nes singuli\`eres. L'existence de supports bien d\'efinis pour les cha\^\i nes
singuli\`eres permet d'introduire une notion de point fixe pour les morphismes de cha\^\i
nes telle que, pour toute fonction continue $f:X\to X$, un point de $X$ soit un point
fixe de $f$ si, et seulement si, c'est un point fixe du morphisme de cha\^\i nes induit
par $f$. Cette approche a \'et\'e utilis\'ee dans [5] pour construire l'indice de point
fixe des fonctions compactes admissibles ; nous l'exposerons \`a la section 3. A la
section 4, nous d\'efinirons une notion d'approximation des fonctions multivoques
compactes semi-continues sup\'erieurement par des morphismes de cha\^\i nes, et la
th\'eorie des points fixes des morphismes de cha\^\i nes permettra d'obtenir des
th\'eor\`emes de point fixe pour ces fonctions multivoques.

Le passage du cas m\'etrisable au cas g\'en\'eral utilise deux ingr\'edients. Le premier
est un foncteur associant \`a tout compact $X$ un espace uniform\'ement contractile
$U(X)$ le contenant et ayant la propri\'et\'e que toute fonction continue de $X$ dans un
espace UC se prolonge contin\^ument \`a $U(X)$. Cette propri\'et\'e permet de r\'eduire
la d\'emonstration de tout th\'eor\`eme de point fixe pour les fonctions compactes
(univoques ou multivoques) des espaces UC au cas particulier des fonctions envoyant
$U(X)$ dans $X$, $X$ compact. Si le th\'eor\`eme qui nous int\'eresse est d\'ej\`a
d\'emontr\'e por les espaces UC m\'etrisables, la v\'erification de ce cas particulier
est simple lorsque le compact $X$ est m\'etrisable. Le deuxi\`eme ingr\'edient, qui
permet de passer des compacts m\'etrisables aux compacts arbitraires dans le cas
particulier, est un type sp\'ecial de syst\`emes projectifs de compacts m\'etrisables.
Ces deux ingr\'edients et la fa\c con de les utiliser seront d\'ecrits aux sections 6 et
7. Les d\'emonstrations qui suiven t ce plan sont parfois longues, mais m\'ecaniques.
Nous donnerons un exemple d'une telle d\'emonstration dans la derni\`ere section.

Si $\U$ est un recouvrement ouvert d'un espace $X$ et $A$ un sous-ensemble de $X$, nous
posons $\St(A,\U)=\bigcup\{U\in\U\,|\,U\cap A\ne\vi\}$, et nous d\'efinissons le
recouvrement ouvert $\St(\U)$ par $\St(\U)=\{\St(U,\U)\,|\,U\in\U\}$. Le recouvrement de
$A$ par les ensembles $A\cap U$ avec $U\in\U$ sera not\'e $\U|A$. Si $\U_1$ et $\U_2$
sont des recouvrements ouverts de $X_1$ et $X_2$ respectivement, nous notons
$\U_1\ti\U_2$ le recouvrement de $X_1\ti X_2$ form\'e par les ensembles $U_1\ti U_2$,
o\`u $U_1$ appartient \`a $\U_1$ et $U_2$ \`a $\U_2$.

Npous notons $|K|$ la r\'ealisation g\'eom\'etrique du complexe simplicial $K$. Si $\si$,
$\si'$ sont deux simplexes de $K$, la notation $\si\le\si'$ signifie que $\si$ est une
face de $\si'$ et $\si<\si'$ que $\si$ est une face propre de $\si'$. Si les sommets
$v_0,\dots,v_k$ de $K$ engendrent un simplexe, nous notons $[v_0,\dots,v_k]$ ce simplexe
et $|v_0,\dots,v_k|$ sa r\'ealisation g\'eom\'etrique. La subdivision barycentrique de
$K$ est not\'ee $K'$. Pour tout simplexe $\si$ de $K$, nous notons $b_\si$ son barycentre
et $\Tr\si$ le sous-complexe de $K'$ form\'e des simplexes $[b_{\si_0},\dots,b_{\si_k}]$
tels que $\si\le\si_0<\dots<\si_k$.

\section{Les r\'etractes absolus de voisinage alg\'ebriques}

Soit $R$ un anneau unitaire. Pour tout espace topologique $X$, nous notons $S(X,R)$ le
complexe des cha\^\i nes singuli\`eres de $X$ \`a coefficients $R$, $H(X,R)$ son groupe
gradu\'e d'homologie et $\wt H(X,R)$ son homologie r\'eduite. Pour tout sous-ensemble $A$
de $X$, nous identifions $S(A,R)$ \`a un sous-complexe de $S(X,R)$. Si $f:X\to Y$ est une
fonction continue, nous notons $f_\#:S(X,R)\to S(Y,R)$ et $f_*:H(X,R)\to H(Y,R)$ les
morphismes qu'elle induit. Nous identifions chaque simplexe singulier $\si$ de $X$ \`a
l'\'el\'ement $1.\si$ de $S(X,R)$, et notons $||c||$ le support d'une cha\^\i ne $c\in
S(X,R)$. Si $\U$ est un recouvrement ouvert de $X$, nous notons $S(X,\U,R)$ le
sous-complexe de $S(X,R)$ engendr\'e par les simplexes singuliers dont l'image est
contenue dans un \'el\'ement de $\U$. Il existe (voir [11]) un morphisme de cha\^\i nes
$Sd:S(X,R)\to S(X,\U,R)$, appel\'e un op\'erateur de subdivision, \'egal \`a l'identit\'e
sur $S(X,\U,R)$ et v\'erifiant $||Sd(c)||\su||c||$ pour toute cha\^\i ne $c\in S(X,R)$ ;
en outre, il existe une homotopie $h$ entre l'identit\'e de $S(X,R)$ et $Sd$ telle que
$||h(c)||\su||c||$ pour toute cha\^\i ne $c$. En particulier, l'inclusion de $S(X,\U,R)$
dans $S(X,R)$ induit un isomorphisme sur l'homologie, par lequel nous identifions
l'homologie de $S(X,\U,R)$ \`a $H(X,R)$ et, si $\ph:S(X,\U,R)\to S(X,R)$ est un morphisme
de cha\^\i nes, nous regardons l'homomorphisme $\ph_*$ qu'il induit sur l'homologie comme
un endomorphisme de $H(X,R)$.

Pour tout complexe simplicial $K$, nous notons $C(K,R)$ le complexe des cha\^\i nes
orient\'ees de $K$ \`a coefficients $R$ (voir [11]), $C_q(K,R)$ le groupe des $q$-cha\^\i
nes de ce complexe et $H(K,R)$ son groupe gradu\'e d'homologie. Pour tout $q\ge0$, nous
identifions $C_q(K,R)$ au $R$-module libre engendr\'e par les $q$-simplexes de $K$ en
fixant un g\'en\'erateur $\si$ de $C_{\dim\si}(\si,R)$ pour chaque simplexe $\si$ de $K$.
Le symbole $\si$ d\'esignera donc \`a la fois un simplexe de $K$ et l'\'el\'ement
correspondant de $C(K,R)$.

Tous les complexes de cha\^\i nes utilis\'es dans cet article sont naturellement
augment\'es, et {\it tous les morphismes de cha\^\i nes entre deux tels complexes seront
suppos\'es pr\'eserver l'augmentation}.

La d\'efinition des r\'etractes absolus de voisinage alg\'ebriques est un analogue de la
caract\'erisation de Lefschetz des r\'etractes absolus de voisinage ([10], th\'eor\`eme
4.1, p 122) dans laquelle les r\'ealisations des complexes simpliciaux relativement \`a
un recouvrement ouvert, qui sont des fonctions, sont remplac\'ees par des r\'ealisations
alg\'ebriques, qui sont des morphismes de cha\^\i nes.

\begin{defi}Soient $X$ un espace topologique, $\U$ un recouvrement ouvert de $X$ et $K$
un complexe simplicial. Une r\'ealisation alg\'ebrique partielle de $K$ relativement \`a
$\U$ est la donn\'ee d'un sous-complexe $L$ de $K$ contenant tous les sommets de $K$ et
d'un morphisme de cha\^\i nes $\mu:C(L,R)\to S(X,R)$ tel que, pour tout simplexe $\si$ de
$K$, il existe un \'el\'ement de $\U$ contenant $||\mu(\tau)||$ pour toute face $\tau$ de
$\si$ appartenant \`a $L$. Si $L=K$, la r\'ealisation alg\'ebrique est dite compl\`ete.
\end{defi}

\begin{defi}Un espace m\'etrisable $X$ est appel\'e un r\'etracte absolu de voisinage
alg\'ebrique, ou RAV alg\'ebrique, si, pour tout recouvrement ouvert $\U$ de $X$, il
existe un recouvrement ouvert $\V$ de $X$ qui est plus fin que $\U$ et tel que, pour tout
complexe simplicial $K$, toute r\'ealisation alg\'ebrique partielle de $K$ relativement
\`a $\V$ se prolonge en une r\'ealisation alg\'ebrique compl\`ete de $K$ relativement \`a
$\U$.
\end{defi}

La d\'efinition d'un RAV alg\'ebrique d\'epend de l'anneau $R$ de coefficients. Lorsque
nous voudrons pr\'eciser cet anneau, nous parlerons d'un $R$-RAV alg\'ebrique.

Pour l'\'etude des RAV alg\'ebriques, il est avantageux d'utiliser aussi certains
CW-complexes particuliers. Par une cellule d'un CW-complexe $M$, nous entendrons toujours
un sous-ensemble ferm\'e de $M$, et nous noterons $\oc{\si}$ l'int\'erieur formel de la
cellule $\si$. Si $M$ est un CW-complexe, nous notons $C(M,R)$ le complexe des
$R$-cha\^\i nes cellulaires de $M$. C'est un module libre ayant une base en
correspondance biunivoque avec les cellules de $M$, le g\'en\'erateur correspondant \`a
une cellule $\si$ de $M$ appartenant \`a $H_{\dim\si}(\si,\dot\si,R)$ ; nous supposerons
toujours avoir fix\'e une telle base, et noterons le g\'en\'erateur correspondant \`a la
cellule $\si$ par la m\^eme lettre. Nous dirons qu'un CW-complexe $M$ est sp\'ecial si
toutes ses cellules sont des sous-complexes et s'il admet une subdivision simpliciale.
Les notions de sommet et de face d'une cellule ont un sens pour tout CW-complexe
sp\'ecial, donc la d\'efinition d'une r\'ealisation alg\'ebrique (partielle) relativement
\`a un recouvrement ouvert s'applique aussi aux CW-complexes sp\'eciaux. Il est prouv\'e
au lemme 1 de [4] que si le recouvrement ouvert $\V$ v\'erifie la condition de la
d\'efinition 2 relativement au recouvrement ouvert $\U$ et \`a tout complexe simplicial,
alors il v\'erifie aussi cette condition relativement au recouvrement $\St(\U)$ et \`a
tout CW-complexe sp\'ecial. Il est donc possible, pour d\'efinir les RAV alg\'ebriques,
d'utiliser aussi bien les CW-complexes sp\'eciaux que les complexes simpliciaux. Les
CW-complexes sp\'eciaux se rencontrent naturellement de deux fa\c cons:

\noindent I) Si $M$ est un CW-complexe sp\'ecial, alors $M\ti I$ est un CW-complexe
sp\'ecial de cellules $\si_0=\si\ti\{0\}$, $\si_1=\si\ti\{1\}$ et $\si_I=\si\ti I$, o\`u
$\si$ parcourt les cellules de $M$. Nous choisissons les g\'en\'erateurs $\si_0$, $\si_1$
et $\si_I$ correspondant \`a ces cellules de fa\c con que la projection envoie $\si_0$ et
$\si_1$ sur $\si$, puis, inductivement sur la dimension de $\si$, de fa\c con que
$\pa\si_I=\si_1-\si_0-(\pa\si)\ti I$ (o\`u, si $c=\sum\ep_i\si^i$, alors $c\ti I=\sum
\ep_i\si_I^i)$. Un morphisme de cha\^\i nes $\ph:C(M\ti I,R)\to C'$ correspond alors \`a
deux morphismes de cha\^\i nes $\ph_0$, $\ph_1$ de $C(M,R)$ dans $C'$ et \`a une
homotopie $h$ entre $\ph_0$ et $\ph_1$ ($\ph_1-\ph_0=\pa h+h\pa)$ d\'efinis par
$\ph_0(\si)=\ph(\si_0)$, $\ph_1(\si)=\ph(\si_1)$ et $h(\si)=\ph(\si_I)$. La d\'efinition
2, appliqu\'ee \`a $M\ti I$, permet donc de construire une homotopie entre deux
morphismes de cha\^\i nes ou de prolonger une homotopie alg\'ebrique d\'efinie sur un
sous-complexe (un analogue alg\'ebrique de la propri\'et\'e classique d'extension des
homotopies).

\noindent II) Pour tout espace topologique $Y$, soit $\Si(Y)$ la r\'ealisation
g\'eom\'etrique (au sens de Giever [8]) du complexe singulier de $Y$. C'est un
CW-complexe sp\'ecial dont les cellules sont en correspondance biunivoque avec les
simplexes singuliers de $Y$ et dont le complexe cellulaire $C(\Si(Y),R)$ est
naturellement isomorphe au complexe singulier $S(Y,R)$. Cela permet d'appliquer la
d\'efinition 2 pour construire des morphismes de cha\^\i nes sur des sous-complexes de
$S(Y,R)$.

De m\^eme que les r\'etractes absolus de voisinage sont d\'efinis en termes de
r\'etractions, les RAV alg\'ebriques peuvent \^etre caract\'eris\'es en termes de
r\'etractions alg\'ebriques.

\begin{defi}Un sous-espace $A$ d'un espace $X$ est appel\'e un r\'etracte de voisinage
alg\'ebrique de $X$ s'il esiste un voisinage $U$ de $A$ dans $X$, un recouvrement ouvert
$\U$ de $U$ et un morphisme de cha\^\i nes $\mu:S(U,\U,R)\to S(A,R)$ v\'erifiant:

\begin{enumerate}
\item[(i)]$\mu(c)=c$ pour $c\in S(A,R)\cap S(U,\U,R)$,

\item[(ii)]pour tout $x\in A$ et tout voisinage $V$ de $x$ dans $A$, il existe un
voisinage $W$ de $x$ dans $X$ tel que $\mu\big(S(W,R)\cap S(U,\U,R)\big)\su
S(V,R)$\footnote{Cette d\'efinition d\'epend aussi de l'anneau $R$. Nous utiliserons
toujours un seul anneau \`a la fois, le m\^eme pour les RAV alg\'ebriques et pour les
r\'etractions alg\'ebriques.}.
\end{enumerate} \end{defi}

Tout r\'etracte de voisinage est un r\'etracte de voisinage alg\'ebrique car le morphisme
$r_\#:S(U,R)\to S(A,R)$ induit par une r\'etraction $r:U\to A$ v\'erifie les conditions
de la d\'efinition pr\'ec\'edente.

\begin{theo}Pour un espace m\'etrisable $X$, les deux conditions suivantes sont
\'equivalentes.

\begin{enumerate}
\item[(i)]$X$ est un RAV alg\'ebrique.

\item[(ii)]$X$ est un r\'etracte de voisinage alg\'ebrique de tout espace m\'etrisable le
contenant comme ferm\'e.
\end{enumerate} \end{theo}

Il r\'esulte du th\'eor\`eme 1 et de la remarque qui le pr\'ec\`ede que tout r\'etracte
absolu de voisinage est un RAV alg\'ebri\-que, mais la classe des RAV alg\'ebri\-ques est
beaucoup plus large que celle des r\'etractes absolus de voisinage, comme le montre en
particulie le th\'eor\`eme 3 ci-dessous.

Le r\'esultat suivant est un analogue de la caract\'erisation des r\'etractes absolus de
voisinage en termes de petites dominations par des complexes simpliciaux ([10],
th\'eor\`eme 6.3, p 139). Il joue un r\^ole essentiel dans notre approche de l'indice de
point fixe.

\begin{theo}Pour un espace m\'etrisable $X$, les conditions suivantes sont
\'equivalentes.

\begin{enumerate}
\item[(i)]$X$ est un RAV alg\'ebrique.

\item[(ii)]Pour tout recouvrement ouvert $\U$ de $X$, il existe un recouvrement ouvert
$\V$ de $X$ plus fin que $\U$, un complexe simplicial $K$, des morphismes de cha\^\i nes
$\psi:S(X,\V,R)\to C(K',R)$, $\zeta:C(K',R)\to S(X,R)$ et une homotopie $h:S(X,\V,R)\to
S(X,R)$ entre l'inclusion et $\zeta\circ\psi$ v\'erifiant

\begin{enumerate}
\item[(a)]Pour tout $V\in\V$, il existe un simplexe $s$ de $K$ tel que $\psi(S(V,R))$
soit contenu dans $C(\Tr s,R)$ $[$et il y a un \'el\'ement de $\U$ contenant $V$ et
$||\zeta(t)||$ pour tout $t\in\Tr s]$.

\item[(b)]$\zeta$ est une r\'ealisation alg\'ebrique compl\`ete de $K'$ relativement \`a
$\U$.

\item[(c)]Pour tout simplexe singulier $\si$ appartenant \`a $S(X,\V,R)$, il existe $U\in
\U$ contenant $||\tau||\cup||\zeta\circ\psi(\tau)||\cup||h(\tau)||$ pour toute face
$\tau$ de $\si$. \end{enumerate}

En outre

\begin{enumerate}
\item[(d)]Si $C$ est un compact fix\'e de $X$, ces objets peuvent \^etre choisis de fa\c
con qu'il existe un voisinage $O$ de $C$ dans $X$ et un sous-complexe fini $L$ de $K'$
tels que $\psi\big(S(X,\V,R)\cap S(O,R)\big)\su C(L,R)$.
\end{enumerate} \end{enumerate} \end{theo}

Ce th\'eor\`eme est prouv\'e dans [4] sans la deuxi\`eme partie (entre crochets) de la
condition (a), qui n'est pas n\'ec\'essaire pour la d\'emonstration de l'implication (ii)
$\Rightarrow$(i). Le fait que la d\'emonstration donn\'ee dans [4] permet aussi d'obtenir
la deuxi\`eme partie de (a) est remarqu\'e dans [5].

Le corollaire suivant permet, en particulier de d\'efinir le nombre de Lefschetz d'une
fonction continue compacte.

\medskip \noindent
{\bf Corollaire.} {\it Soient $R$ un anneau n\oe th\'erien et $X$ un $R$-RAV
alg\'ebrique. Pour tout compact $C$ de $X$, il existe un voisinage $O$ de $C$ dans $X$
tel que l'image de l'homomorphisme de $H(O,R)$ dans $H(X,R)$ induit par l'inclusion soit
engendr\'ee par un nombre fini d'\'el\'ements.}

\begin{theo}Tout espace m\'etrisable ULC est un RAV alg\'ebrique. \end{theo}

La d\'emonstration de ce th\'eor\`eme est longue et d\'elicate. Nous l'esquisserons dans
le cas d'un sous-ensemble convexe compact d'un e.v.t. ; ce cas particulier contient les
id\'ees les plus importantes et les plus originales de la d\'emonstration.

Posons $\Dt_n=\{(t_0,\dots,t_n)\in\Bbb R^{n+1}\,|\,t_i\ge0 \text{ pour tout $i$ et
}\sum_{i=0}^nt_i=1\}$, et, pour tout $n\ge1$, soit $\la_n:X^{n+1}\ti\Dt_n\to X$ la
fonction continue d\'efinie par
$$
\la_n(x_0,\dots,x_n;t_0,\dots,t_n)=\sum_{i=0}^nt_ix_i.
$$

Plongeons $X$ dans le cube de Hilbert $Y=\prod_{i=1}^\ii I_i$, o\`u $I_i=[0,1]$ pour tout
$i$, muni de la distance d\'efinie, pour $x=(x_i)$ et $y=(y_i)$ par
$d(x,y)=\sum_{i=1}^\ii 2^{-i}|x_i-y_i|$. Pour $x\in X$ et $\ep>0$, $B(x,\ep)$ d\'esignera
la boule ouverte de centre $x$ et de rayon $\ep$ dans $X$ pour la distance induite par
$d$. Pour $n\ge1$, soient $Y^n=\prod_{i=1}^nI_i$ et $p_n$ la projection de $Y$ sur $Y^n$.

Nous voulons montrer que $X$ v\'erifie la condition (ii) du th\'eor\`eme 2. La t\^ache
est simplifi\'ee par le lemme 4 de [4] gr\^ace auquel nous n'avons pas \`a nous soucier
de construire l'homotopie $h$. Par compacit\'e, il suffit de consid\'erer le cas o\`u,
pour $\ep>0$ fix\'e, $\U$ est le recouvrement form\'e par les boules $B(x,\ep)$, $x\in
X$. Puisque $\la_2(x,x,x;t_0,t_1,t_2)=x$ pour tout $x\in X$ et tout
$(t_0,t_1,t_2)\in\Dt_2$, nous pouvons trouver $0<\ep_1<\ep$ tel que
\begin{equation}
\diam\la_2(B(x,\ep_1)^3\ti\Dt_2)<\ep/4 \text{ pour tout }x\in X.
\end{equation}

Soit $n>0$ tel que $2^{-n}<\ep_1/8$. Il existe $\dt>0$ v\'erifiant
\begin{equation}
\diam\la_{n+1}(B(x,\dt)^{n+2}\ti\Dt_{n+1})<\ep_1/2 \text{ pour tout }x\in X.
\end{equation}

Soit $k>n$ tel que $2^{n+1-k}<\dt$. D\'ecomposons chaque intervalle $I_i$ en $2^k$
sous-intervalles ferm\'es de longueur $2^{-k}$, et munissons $Y^n$ et $Y^k$ des
d\'ecompositions cellulaires produits. Soit $N$ le sous-complexe de $Y^n$ form\'e des
cellules $\si$ telles que $p_n^{-1}(\oc{\si})\cap X\ne\vi$ et de toutes les faces de
telles cellules, et soit $M$ le sous-complexe de $Y^k$ form\'e des cellules $\tau$ telles
que $p_k^{-1}(\oc{\tau})\cap X\ne\vi$ et de toutes les faces de telles cellules. La
restriction \`a $M$ de la projection de $Y^k$ sur $Y^n$ est une application cellulaire
que nous notons $\pi$. Si les points $(x_i)$ et $(y_i)$ appartiennent \`a une cellule
$\si$ de $N$, alors $|x_i-y_i|\le2^{-k}$ pour $i\le n$, d'o\`u
\begin{equation}
\diam p_n^{-1}(\si)<2^{-n}+2^{-k}<2.2^{-n}\text{ pour toute cellule $\si$ de $N$.}
\end{equation}
De m\^eme,
\begin{equation}
\diam p_k^{-1}(\tau)\le2.2^{-k}\text{ pour toute cellule $\tau$ de $M$.}
\end{equation}
Pour toute cellule $\si$ de $N$, soit $M(\si)$ le sous-complexe de $M$ form\'e des
cellules $\tau$ telles que $\pi(\tau)=\si$ et de toutes les faces de telles cellules.
$M(\si)$ est de la forme $M(\si)=\si\ti L(\si)$, o\`u $L(\si)$ est un sous-complexe de la
d\'ecomposition CW produit de $\prod_{i=n+1}^kI_i$. Soient $L_j$, $j\in J_\si$, les
composantes de $L(\si)$; alors les $M_j=\si\ti L_j$ sont les composantes de $M(\si)$.
Soit $J$ la r\'eunion disjointe des $J_\si$, o\`u $\si$ parcourt les cellules de $N$. Si
$\si'$ est une face de $\si$, alors $M(\si)\cap\pi^{-1}(\si')$ est contenu dans
$M(\si')$, de sorte que, pour tout $j\in J_\si$, il existe un unique $j'\in J_{\si'}$ tel
que $M_j\cap\pi^{-1}(\si')=\si'\ti L_j\su M_{j'}$. Nous \'ecrivons $j'\le j$ si $j'\in
J_{\si'}$ et $j\in J_\si$ avec $\si'\le\si$ et $M_j\cap\pi^{-1}(\si')\su M_{j'}$.

L'ensemble des sommets du complexe simplicial $K$ est $J$ et $s=[j_0,\dots,j_m]$ est un
simplexe de $K$ si la num\'erotation peut \^etre choisie de fa\c con que $j_0\le\dots\le
j_m$. Si tous les $M(\si)$ sont connexes, alors $K$ est isomorphe \`a la subdivision
barycentrique de $N$ (qui est un complexe simplicial puisque toutes les cellules de $K$
sont convexes). Dans le cas g\'en\'eral, la fonction $q$ qui fait correspondre \`a chaque
\'el\'ement de $J_\si$ le barycentre de la cellule $\si$ est une application simpliciale
de $K$ dans la subdivision barycentrique de $N$ dont la restriction \`a chaque simplexe
de $K$ est injective, donc $\dim K\le n$.

Pour construire $\V$, nous avons besoin d'un recouvrement auxiliaire $\G=\{G_j\,|\,j\in
J\}$ index\'e par les sommets de $K$. Pour $j\in J_\si$, $p_k^{-1}(M_j)$ est contenu dans
$p_n^{-1}(\si)$, donc de diam\`etre inf\'erieur \`a $2.2^{-n}$ d'apr\`es (3). Pour $j\in
J_\si$, posons $\dim j=\dim\si$, et, pour tout entier $r\ge0$, soit $J_r=\{j\in
J\,|\,\dim j=r\}$. Les ensembles $p_k^{-1}(M_j)\cap X$ avec $j$ dans $J_0$ sont ferm\'es,
disjoints et de diam\`etres inf\'erieurs \`a $2.2^{-n}$, donc nous pouvons trouver une
famille $\{G_j\,|\,j\in J_0\}$ d'ouverts disjoints de $X$ v\'erifiant $p_k^{-1}(M_j)\cap
X\su G_j$ et $\diam G_j<2.2^{-n}$ pour tout $j\in J_0$. Supposant les $G_{j'}$ construits
pour $\dim j'<r$, soit $O_r=\bigcup_{\dim j'<r}G_{j'}$. Alors $\{p_k^{-1}(M_j)\cap(X\sm
O_r)\,|\,j\in J_r\}$ est une famille finie de ferm\'es disjoints de diam\`etres
inf\'erieurs \`a $2.2^{-n}$, donc nous pouvons trouver une famille $\{G_j\,|\,j\in J_r\}$
d'ouverts disjoints v\'erifiant $p_k^{-1}(M_j)\cap(X\sm O_r)\su G_j$ et $\diam
G_j<2.2^{-n}$ pour tout $j\in J_r$. Si les $G_j$ sont construits avec suffisamment de
soin, le nerf de $\G$ est un sous-complexe de $K$. Soit $\R$ l'ensemble des simplexes de
$K$. Pour $s=[j_0,\dots,j_r]\in\R$, posons $V_r=G_{j_0}\cap \dots\cap G_{j_r}$. Alors
$\V=\{V_s\,|\,s\in\R\}$ est le recouvrement cherch\'e de $X$.

Le morphisme $\psi:S(X,\V,R)\to C(K',R)$ est construit de fa\c con que $\psi(S(V_s,R))\su
C(\Tr s,R)$ pour tout $s$. Pour tout simplexe $s$ de $K$, soit $d(s)$ la dimension de
$\Tr s$. Comme le nerf de $\G$ est un sous-complexe de $K$, si $s_1$ et $s_2$ sont deux
simplexes distincts de $K$ tels que $V_{s_1}\cap V_{s_2}\ne\vi$, alors $s_1\cup s_2$ est
un simplexe de $K$ et $V_{s_1\cup s_2}\su V_{s_1}$. Par suite, si $S(V_{s_1},R)\cap
S(V_{s_2},R)\ne\{0\}$, alors $s_1\cup s_2$ est un simplexe de $K$ et nous avons $\Tr
s_1\cup s_2\su\Tr s_1$ et $d(s_1)>d(s_1\cup s_2)$ si $s_2$ n'et pas une face de $s_1$.
L'acyclicit\'e des complexes $\Tr s$ permet alors de construire $\psi$ en d\'efinissant
par r\'ecurrence sur $d(s)$ sa restriction aux complexes $S(V_s,R)$ dont $S(X,\V,R)$ est
la somme.

La construction de $\zeta:C(K',R)\to S(X,R)$ est plus compliqu\'ee. Nous avons besoin de
quelques objets auxiliaires: une d\'ecomposition CW $\wt K$ de $|K|$ \footnote { Bien que
$\wt K=|K|$, nous utilisons cette nouvelle notation pour distinguer les complexes
cellulaires $C(\wt K,R)$ et $C(|K|,R)$.}, un sous-complexe $Q$ de $M$, un CW-complexe $P$
contenant $Q$ et une fonction continue $f:P\to X$. Le morphisme $\zeta$ sera donn\'e par
une composition

$$ C(K',R)\stackrel{\nu}{\longrightarrow}C(\wt K,R)\stackrel{\ta}{\longrightarrow}
C(P,R)\stackrel{\mu}{\longrightarrow}S(P,R)\stackrel{f_*}{\longrightarrow}S(X,R)
$$

Soient $j\in J$ et $\si$ la cellule de $N$ telle que $j\in J_\si$. Pour toute face $\si'$
de $\si$, soit $j(\si')$ l'\'el\'ement de $J_{\si'}$ tel que $M_j\cap\pi^{-1}(\si')\su
M_{j(\si')}$, et soit $c_j$ le sous-complexe plein de $K$ engendr\'e par les $j(\si')$,
$\si'\le\si$. L'application simpliciale $q$ envoie $c_j$ isomorphiquement sur la
subdivision barycentrique de $\si$, donc $|c_j|$ est une cellule, et les $|c_j|$, $j\in
J$, forment une d\'ecomposition CW sp\'eciale de $|K|$, que nous notons $\wt K$. Pour
tout simplexe $s'$ de $K'$, soit $|c_{j(s')}|$ la plus petite cellule de $\wt K$
contenant $|s'|$. Soit $\nu:C(K',R)\to C(\wt K,R)$ un morphisme de cha\^\i nes tel que
$\nu(s')\in C(|c_{j(s')}|,R)$ pour tout simplexe $s'$ de $K'$.

Pour $j\in J$, soit $S_j$ le $1$-squelette de $L_j$; puisque $L_j$ est connexe, $S_j$
aussi. Si $\si$ est la cellule telle que $j\in J_\si$, posons $Q_j=\si\ti S_j$; c'est un
complexe de dimension au plus $\dim\si+1$, donc $Q=\bigcup_{j\in J}Q_j$ est un
sous-complexe de $M$ de dimension au plus $n+1$. Soit $\{v_j\,|\,j\in J\}$ une famille de
points distincts n'appartenant pas \`a $M$. Pour $j\in J$, soit $C_j=\bigcup\{Q_{j'}\,|\,
j'\le j\}$, et soit $D_j$ le c\^one de base $C_j$ et de sommet $v_j$. Pour $j'<j$, soit
$E(j',j)$ le joint de $C_{j'}$ et du $1$-simplexe g\'eom\'etrique $|v_{j'}v_j|$. Posons
$$
P=Q\cup\bigcup_{j\in J}D_j\cup\bigcup_{j'<j}E(j',j).
$$

C'est un CW-complexe sp\'ecial dont les cellules sont celles de $P$, les c\^ones de
sommet $v_j$ sur les cellules de $C_j$ et les joints de $|v_{j'}v_j|$ avec les cellules
de $C_{j'}$. Soit $\mu:C(P,R)\to S(P,R)$ un morphisme de cha\^\i nes tel que $\mu(\tau)
\in S(\tau,R)$ pour toute cellule $\tau$ de $P$.

Les cellules de $M$ \'etant convexes, sa subdivision barycentrique en est une
triangulation. Nous notons $b_\tau$ le barycentre de la cellule $\tau$ de $M$. Pour
chaque cellule $\tau$ de $M$, il existe une cellule $\tau'$ de $M$ dont $\tau$ est face
et telle que $p_k^{-1}(\oc{\tau'})\cap X\ne\vi$; fixons une telle cellule $\tau'$ et un
point $x(\tau)$ de $p_k^{-1}(\oc{\tau'})\cap X$. Soit $t=|b_{\tau_0},\dots,b_{\tau_r}|$
un simplexe de la subdivision barycentrique de $Q$. Pour un point
$x=\sum_{i=0}^ra_ib_{\tau_i}$ de $t$ ($(a_0,\dots,a_r)\in\Dt_r)$, posons
$f(x)=\sum_{i=0}^ra_ix(\tau_i)$. Pour tout $j\in J$, choisissons une cellule $\tau(j)$ de
$C_j$ et posons $f(v_j)=x(\tau(j))$. Tout point $y\in D_j$ (resp. $z\in E(j',j)$) peut
s'\'ecrire $y=t_0x+t_1v_j$ (resp. $z=t_0x+t_1v_j+t_2v_{j'}$) o\`u $x$ appartient \`a
$C_j$ (resp. $C_{j'}$) et $(t_0,t_1)$ \`a $\Dt_1$ (resp. $(t_0,t_1,t_2)$ \`a $\Dt_2$).
Posons $f(y)=t_0f(x)+t_1f(v_j)$ et $f(z)=t_0f(x)+t_1f(v_j)+t_2f(v_{j'})$.

Pour tout $j\in J$, posons $\Pi_j=\bigcup_{j'\le j}D_{j'}\cup\bigcup_{j''<j'\le
j}E(j'',j')$. C'est un sous-complexe de $P$ contenant $Q_j$ et $j'\le j$ entra\^\i ne
$\Pi_{j'}\su\Pi_j$. Pour tout $j$, nous avons

\begin{equation}
\diam f(\Pi_j)<\ep/2.
\end{equation}

En effet, pour $j''<j'\le j$, les intersections $D_j\cap D_{j'}$ et $D_j\cap E(j'',j')$
sont non vides, donc il suffit de v\'erifier que chacun des ensembles $f(D_j)$ et
$f(E(j',j))$ a un diam\`etre inf\'erieur \`a $\ep/4$. D'apr\`es (1), il suffit pour cela
de constater que les points $f(v_{j'})$, $j'\le j$, et $f(x)$, $x\in C_j$, sont dans la
boule $B(f(v_j),\ep_1)$. Soit $\si$ tel que $j\in J_\si$. Alors $C_j$ est contenu dans
$\pi^{-1}(\si)$, donc $\tau'\cap\pi^{-1}(\si)\ne\vi$ pour toute cellule $\tau$ de $C_j$.
Utilisant (3) et (4), nous constatons que, si $\tau_1$ et $\tau_2$ sont deux cellules de
$C_j$, alors
$$
\begin{aligned}
d(x(\tau_1),x(\tau_2)) &\le\diam p_k^{-1}(\tau'_1)+\diam p_n^{-1}(\si)+\diam p_k^{-1}(\tau'_2) \\
 &<2.2^{-n}+4.2^{-k}<4.2^{-n}<\ep_1/2.
 \end{aligned}$$

Comme $\tau(j')$ appartient \`a $C_{j'}\su C_j$ pour $j'\le j$, nous avons en particulier
$d(f(v_j),f(v_{j'}))<\ep_1/2$. Si $t=|b_{\tau_0},\dots,b_{\tau_r}|$,
$\tau_0<\dots<\tau_r$, est un simplexe de la subdivision barycentrique de $C_j$, nous
avons $\tau_0\su\tau'_0\cap\tau'_i$ pour $0\le i\le r$, donc $d(x(\tau_0),x(\tau_i))<
4.2^{-k}<\dt$ d'apr\`es (4). Puisque $\dim Q\le n+1$, nous avons $r\le n+1$, et (2)
entra\^\i ne que le diam\`etre de $f(t)$ est inf\'erieur \`a $\ep_1/2$; comme
$f(b_{\tau_0})=x(\tau_0)$, nous obtenons, pour $x\in t$,
$$
d(f(v_j),f(x))\le d(f(v_j),x(\tau_0))+d(x(\tau_0),f(x))<\ep_1.
$$

Notant $c_j$ le g\'en\'erateur de $C(\wt K,R)$ correspondant \`a la cellule $|c_j|$, le
morphisme $\ta$ est construit de fa\c con que $\ta(c_j)$ appartienne \`a $C(\Pi_j,R)$. Il
ne reste plus alors qu'\`a v\'erifier la condition (b) et la version de la condition (c)
dans laquelle l'homotopie $h$ est omise. Si $s'$ est un simplexe de $K'$ et $|c_{j(s')}|$
la plus petite cellule de $\wt K$ contenant $|s'|$, les constructions pr\'ec\'edntes
garantissent que $\zeta(s'')$ appartient \`a $S(f(\Pi_{j(s')}),R)$ pour toute face $s''$
de $s'$; d'apr\`es (5), $f(\Pi_{j(s')})$ est contenu dans une boule $B(x,\ep)$, d'o\`u
(b). Soit $V_s=G_{j_0}\cap\dots\cap G_{j_r}$ un \'el\'ement non vide de $\V$, o\`u
$s=[j_0,\dots,j_r]$ ($j_0<\dots<j_r$) est un simplexe de $K$, et soit $\si$ tel que
$j_r\in J_\si$. Fixons un point $x_s$ de $G_{j_r}\cap p_k^{-1}(M_{j_r})$. Si $s'$ est un
simplexe de $\Tr s$, alors $|c_{j(s')}|$ contient $|c_{j_r}|$ et $j_r\le j(s')$, donc
$\zeta(C(\Tr s,R))$ est contenu dans $S(f(\bigcup_{j_r\le j}\Pi_j),R)$ et, pour
v\'erifier la version simplifi\'ee de (c), il suffit de montrer que $B(x_s,\ep)$ contient
$V_s\cup f(\bigcup_{j_r\le j}\Pi_j)$. Comme $\diam G_{j_r}<2.2^{-n}<\ep$, $V_s$ est
contenu dans $B(x_s,\ep)$. Soient $\tau$ une cellule de $C_{j_r}$ et $\tau'$ telle que
$x(\tau)\in p_k^{-1}(\tau')$. Comme $p_n^{-1}(\si)$ contient $x_s$ ainsi que
$p_k^{-1}(\tau)\su p_k^{-1}(\tau')$, nous avons
$$
d(x_s,x(\tau))\le\diam p_n^{-1}(\si)+\diam p_k^{-1}(\tau')<2.2^{-n}+2.2^{-k}<3.2^{-n}.
$$

Pour $j\ge j_r$, $\Pi_j$ contient $\Pi_{j_r}$ et, comme $x(\tau)$ appartient \`a
$f(\Pi_{j_r})$, nous obtenons, pour $y\in \Pi_j$,
$$
\begin{aligned}
d(x_s,f(y))&\le d(x_s,f(\Pi_{j_r}))+\diam f(\Pi_{j_r}) \\
 &<d(x_s,x(\tau))+\ep/2<3.2^{-n}+\ep/2 <\ep, \end{aligned}
$$
d'o\`u l'inclusion $f(\bigcup_{j_r\le j}\Pi_j)\su B(x_s,\ep)$.

La construction de $\ta$ est le seul point de cette d\'emonstration o\`u il est
indispensable d'utiliser des morphismes de cha\^\i nes plut\^ot que des fonctions
continues. En effet, si $\{H_s\,|\,s\in\R\}$ est un recouvrement ferm\'e de $X$ tel que
$H_s\su V_s$ pour tout $s$, on peut trouver une fonction continue $g:X\to|K|$ telle que
$g(H_s)\su|\Tr s|$ pour tout $s$. S'il \'etait toujours possible de construire une
fonction continue $\al:\wt K\to P$ telle que $\al(|c_j|)\su\Pi_j$ pour tout $j$, on en
d\'eduirait que $X$ est un r\'etracte absolu de voisinage. Il est instructif de comparer
les deux cas. Nous sommes dans la situation suivante: \`a tout $j\in J$ est associ\'e un
sous-complexe $Q_j$ de $Q$ de fa\c con que $Q_{j'}\su Q_j$ pour $j'\le j$. Les complexes
$Q_j$ sont connexes, mais leurs groupes d'homologie et d'homotopie sont a priori
arbitraires. Nous voulons construire un CW-complexe $P$ et des sous-complexes $\Pi_j$ de
$P$ v\'erifiant $Q_j\su\Pi_j$ et $\Pi_{j'}\su\Pi_j$ pour $j'\le j$, de fa\c con qu'il
existe soit une fonction contine $\al:\wt K\to P$ v\'erifiant $\al(|c_j|)\su\Pi_j$ pour
tout $j$, soit un morphisme de cha\^\i nes $\ta:C(\wt K,R)\to C(P,R)$ v\'erifiant
$\ta(c_j)\in C(\Pi_j,R)$ pour tout $j$. Le reste de la d\'emonstration est ensuite
adapt\'e \`a la nature de $P$.

Examinons d'abord le cas des fonctions. Si $\dim j\,( =\dim c_j)=0$, envoyons $|c_j|$ sur
un point de $Q_j$. Si $\dim j=1$, la connexit\'e du complexe $Q_j$ garantit que $\al$ se
prolonge en une fonction de $|c_j|$ dans $Q_j$. Nous pouvons donc prendre $\Pi_j=Q_j$
pour $\dim j\le 1$. Si $\dim j=2$, $\al$ est d\'efinie sur le bord de $|c_j|$ et \`a
valeurs dans $Q_j$, mais ne se prolonge pas n\'ec\'essairement en une fonction continue
de $|c_j|$ dans $Q_j$; pour pouvoir la prolonger, nous ajoutons \`a $Q_j$ des c\^ones sur
des sous-ensembles de $Q_j$, ce qui d\'efinit le complexe $\Pi_j$. Si $\dim j=3$, $\al$
est d\'efinie sur le bord de $|c_j|$ et \`a valeurs dans
$Z_j=Q_j\cup\bigcup_{j'<j}\Pi_{j'}$ et, pour pouvoir la prolonger, nous construisons
$\Pi_j$ en ajoutant \`a $Z_j$ des c\^ones sur certains sous-ensembles de $Z_j$; cette
op\'eration introduit des c\^ones sur des sous-ensembles qui sont eux-m\^emes des
c\^ones, c'est \`a dire des joints de sous-ensembles de $Q_j$ avec des $1$-simplexes.
Cette op\'eration se r\'ep\`ete \`a chaque \'etape: pour construire la restriction de
$\al$ \`a une cellule $|c_j|$ de dimension $r$, nous devons ajouter \`a $Q_j$ des joints
de sous-ensembles de $Q_j$ avec des $(r-1)$-simplexes. La dimension maximale des
simplexes avec lesquels nous formons des joints d\'epend de $n$, donc de $\ep$, et il est
impossible de la fixer une fois pour toutes au d\'ebut de la d\'emonstration (ce qui est
en fait la cl\'e de cette d\'emonstration).

Consid\'erons maintenant le cas des morphismes de cha\^\i nes. Si $\dim j=0$, prenons
pour $\ta(c_j)$ le g\'en\'erateur correspondant \`a un sommet de $Q_j$. Si $\dim  j=1$,
$Q_j$ contient une $1$-cha\^\i ne $\ta(c_j)$ telle que $\pa\ta(c_j)=\ta(\pa c_j)$. Si
$\dim j=2$, nous construisons encore $\Pi_j$ en ajoutant \`a $Q_j$ des c\^ones sur
certains de ses sous-ensembles. La cha\^\i ne $\ta(c_j)$ fait alors intervenir des
cha\^\i nes coniques $vd$, joint d'un sommet $v$ et d'une cha\^\i ne $d\in C(Q_j,R)$. Si
$\dim j=3$, $\Pi_j$ s'obtient comme pr\'ec\'edemment en ajoutant \`a $Z_j$ des c\^ones
sur certains de ses sous-ensembles et, dans la d\'efinition de $\ta(c_j)$, apparaissent
des cha\^\i nes $vv'd$, joint d'un $1$-simplexe $vv'$ et d'une cha\^\i ne $d\in
C(Q_j,R)$. C'est ici que la situation se distingue du cas pr\'ec\'edent. En effet, soit
$vd$ une cha\^\i ne conique, et soit $v'$ un sommet distinct de $v$. La cha\^\i ne $vv'd$
a pour bord $v'd-vd+vv'\pa d$. Ajoutant ce bord \`a $vd$, nous obtenons $v'd+vv'\pa d$;
cette op\'eration nous permet donc de changer de sommet en ajoutant un terme correcteur
qui s'annule sur les cycles. Gr\^ace \`a cette simple manipulation alg\'ebrique, il n'est
pas n\'ec\'essaire de prendre des joints avec des simplexes de dimensions sup\'erieures
\`a un, et le complexe $P$ peut \^etre construit comme indiqu\'e plus haut. Nous
renvoyons le lecteur \`a [4] pour la d\'efinition de $\ta(c_j)$ et la v\'erification de
la relation $\ta\pa=\pa\ta$, qui sont assez longues.

En r\'esum\'e, c'est une propri\'et\'e alg\'ebrique \'el\'ementaire, la possibilit\'e de
changer le sommet d'une cha\^\i ne conique en introduisant un terme correcteur nul sur
les cycles, qui permet de prouver que les espaces ULC m\'etrisables sont des RAV
alg\'ebriques. Cette propri\'et\'e n'a aucun analogue homotopique, et il existe des
convexes m\'etrisables qui ne sont pas des r\'etractes absolus, ni m\^eme des extenseurs
absolus pour la classe des compacts m\'etrisables. Les th\'eor\`emes de point fixe pour
les fonctions compactes des r\'etractes absolus de voisinage sont de nature homologique
et, comme nous le verrons ci-dessous, s'\'etendent aux RAV alg\'ebriques. En particulier,
le fait que les convexes m\'etrisables ont la propri\'et\'e du point fixe pour les
fonctions compactes alors que ce ne sont pas n\'ec\'essairement des extenseurs absolus
pour la classe des compacts m\'etrisables est une manifestation en dimension infinie des
diff\'erences de propri\'et\'es entre homologie et homotopie.

\section{L'indice de point fixe pour les RAV alg\'ebriques}

{\it Dans toute cette section, $R$ d\'esigne un anneau principal}.

Soient $X$ un espace topologique s\'epar\'e, $A$ un sous-ensemble de $X$, $\U$ un
recouvrement ouvert de $A$ et $\ph:S(A,\U,R)\to S(X,R)$ un morphisme de cha\^\i nes. Un
point $x$ de $A$ est appel\'e un {\it point fixe} de $\ph$ si, pour tout voisinage $V$ de
$x$ dans $X$, il existe une cha\^\i ne $c\in S(V,R)\cap S(A,\U,R)$ telle que
$||\ph(c)||\cap V\ne\vi$. L'ensemble des points fixes de $\ph$ est not\'e $\Fix\ph$; ce
sous-ensemble est ferm\'e dans $A$. Le morphisme $\ph$ est dit compact s'il existe un
compact $C$ de $X$ telque $S(C,R)$ contienne l'image de $\ph$. Des d\'efinitions
analogues s'appliquent aux homotopies alg\'ebriques $h:S(A,\U,R)\to S(X,R)$.

Si $f:A\to X$ est une fonction continue, alors un \'el\'ement $x$ de $A$ est un point
fixe de $f$ si, et seulement si, c'est un point fixe du morphisme $f_\#:S(A,R)\to
S(X,R)$. Les points fixes des fonctions peuvent donc \^etre \'etudi\'es au niveau des
morphismes induits. Cette approche est particuli\`erement adapt\'ee aux RAV
alg\'ebriques; suivant [5], nous indiquerons ici comment elle permet de construire dans
ces espaces un indice de point fixe pour les fonctions compactes. La premi\`ere t\^ache
est de d\'efinir l'indice de point fixe des morphismes de cha\^\i nes. Nous avons besoin
pour cela de quelques rappels alg\'ebriques.

Un module (gradu\'e ou non) est de type fini s'il est engendr\'e par un nombre fini
d'\'el\'ements. Un morphisme $\ph:C\to C'$ entre $R$-modules (gradu\'es ou non) est de
type fini si son image est de type fini. Le corollaire du th\'eor\`eme 2 entra\^\i ne que
si $\ph:S(Y,\U,R)\to S(X,R)$ est un morphisme de cha\^\i nes compact, o\`u $X$ est un
$R$-RAV alg\'ebrique, alors l'homomorphisme $\ph_*$ induit par $\ph$ sur l'homologie est
de type fini. Si $\ph$ est un endomorphisme de type fini d'un $R$-module $C$, il est
possible de d\'efinir sa trace. Dans le cas o\`u $R$ est un corps, cela est expliqu\'e au
\S 15 de [9]. Si $R$ est un anneau principal de corps des fractions $Q$, la trace de
$\ph$ est d\'efinie comme la trace de l'endomorphisme $\ph\otimes id$ du $Q$-espace
vectoriel $Q\otimes_RC$.

Soit $X$ un $R$-RAV alg\'ebrique. Un morphisme de cha\^\i nes $\ph:S(U,\U,R)\to S(X,R)$,
o\`u $U$ est un ouvert de $X$ et $\U$ un recouvrement ouvert de $U$, sera dit {\it
admissible} s'il est compact et si $\Fix\ph$ est un sous-ensemble compact de $U$. De
m\^eme, une homotopie $h:S(U,\U,R)\to S(X,R)$ sera dite admissible si elle est compacte
et si $\Fix h$ est un sous-ensemble compact de $U$.

L'indice de point fixe $\ind(\ph,U)$ d'un morphisme admissible $\ph:S(U,\U,R)\ \to
S(X,R)$ est d\'efini comme suit. Prenons des voisinages ouverts $A$ et $B$ de $X\sm U$ et
$\Fix\ph$ respectivement tels que $\ov A\cap\ov B=\vi$, puis un recouvrement ouvert $\V$
de $X$ v\'erifiant

\begin{itemize}
\item[($\al$)]$\St(\ov A,\V)\cap\St(\ov B,\V)=\vi$,

\item[($\bt$)]tout \'el\'ement de $\V$ rencontrant $X\sm A$ est contenu dans un
\'el\'ement de $\U$,

\item[($\ga$)]si $V$ est un \'el\'ement de $\V$ tel que $V\sm(A\cup B)\ne\vi$, alors
$||\ph(c)||\cap\St(V,\St(V,\V))=\vi$ pour toute cha\^\i ne $c\in S(V,R)$.
\end{itemize}

La possibilit\'e d'obtenir ($\ga$) r\'esulte du fait que tout point du ferm\'e
$X\sm(A\cup B)$ appartient \`a $U\sm\Fix\ph$, donc a un voisinage $O$ qui est contenu
dans un \'el\'ement de $\U$ et tel que $||\ph(c)||\cap O=\vi$ pour toute $c\in S(O,R)$.

Soit $C$ un  compact de $X$ tel que $S(C,R)$ contienne l'image de $\ph$. Prenons un
recouvrement ouvert $\W$ de $X$, un complexe simplicial $K$, des morphismes de cha\^\i
nes $\psi:S(X,\W,R)\to C(K',R)$, $\zeta:C(K',R)\to S(X,R)$ et une homotopie $h:S(X,\W,R)
\to S(X,R)$ v\'erifiant les conditions (a)-(d) du th\'eor\`eme 2 relativement \`a $\V$ et
$C$. D\'efinissons un endomorphisme $\ta_A$ du groupe gradu\'e $S(X,R)$ en posant, pour
tout simplexe singulier $\si$ de $X$,
$$
\ta_A(\si)=\begin{cases}0 &\text{si }||\si||\cap A\ne\vi \\
\si &\text{sinon}. \end{cases}
$$

L'image de $\ta$ est contenue dans $S(U,R)$, donc si $Sd_1:S(U,R)\to S(U,\U,R)$ et $Sd_2:
S(X,R)\to S(X,\W,R)$ sont des op\'erateurs de subdivision, alors
$$
\xi=\psi\circ Sd_2\circ\ph\circ Sd_1\circ\ta_A\circ\zeta
$$
est un endomorphisme du module gradu\'e $C(K',R)$. Comme $\im\ph\su S(C,R)$, la condition
(d) du th\'eor\`eme 2 garantit que $\xi$ est de type fini, donc $\La(\xi)$ est d\'efini,
et nous posons $\ind(\ph,U)=\La(\xi)$, ce qui est justifi\'e par le lemme suivant.

\begin{lem}(i) $\ind(\ph,U)$ ne d\'epend pas du choix de $A$, $B$, $C$, $\V$, $\W$, $K$,
$\psi$, $\zeta$, $Sd_1$ et $Sd_2$.

(ii) Si $\U'$ est un recouvrement ouvert de $U$ plus fin que $\U$ et si $\ph'$ est la
restriction de $\ph$ \`a $S(U,\U',R)$, alors $\ind(\ph',U)=\ind(\ph,U)$.
\end{lem}
Les propri\'et\'es de cet indice sont analogues aux propri\'et\'es habituelles des
indices de point fixe des fonctions continues.

\begin{theo}L'indice de point fixe pour les morphismes admissibles $\ph:S(U,\U,R)\to
S(X,R)$, o\`u $U$ est un ouvert d'un $R$-RAV alg\'ebrique $X$, a les propri\'et\'es
suivantes:
\begin{itemize}
\item[(I)] {\em (Normalisation)} Si $U=X$, alors $\ind(\ph,X)=\La(\ph_*)$.

\item[(II)] {\em (Additivit\'e)} Si $V_1$, $V_2$ sont deux ouverts disjoints contenus
dans $U$ tels que $\Fix\ph\su V_1\cup V_2$, alors
$\ind(\ph,U)=\ind(\ph,V_1)+\ind(\ph,V_2)$.

\item[(III)] {\em (Existence)} Si $\ind(\ph,U)\ne0$, alors $\Fix\ph\ne\vi$.

\item[(IV)] {\em (Homotopie)} Si $\ph_0$, $\ph_1:S(U,\U,R)\to S(X,R)$ sont admissibles et
s'il existe une homotopie admissible $\mu:S(U,\U,R)\to S(X,R)$ entre $\ph_0$ et $\ph_1$,
alors $\ind(\ph_0,U)=\ind(\ph_1)$.

\item[(V)] {\em (Contraction)} Soit $Y$ un sous-espace de $X$ qui est aussi un $R$-RAV
alg\'ebrique. Supposons qu'il existe un compact $C$ de $Y$ tel que $S(C,R)$ contienne
$\im\ph$. Si $\ph_Y:S(U\cap Y,\U|U\cap Y,R)\to S(Y,R)$ est induit par $\ph$, alors
$\ind(\ph,U)=\ind(\ph_Y,U\cap Y)$.
\end{itemize} \end{theo}

L'indice de point fixe pour les morphismes admissibles a aussi une propri\'et\'e
correspondant \`a la multiplicativit\'e de l'indice de point fixe des fonctions. Soient
$\Phi:S(\,.\,\ti\,.\, ,R)\to S(\,.\, ,R)\otimes S(\,.\, ,R)$ et $\Psi:S(\,.\, ,R)\otimes
S(\,.\, ,R)\to S(\,.\,\ti\,.\, ,R)$ des morphismes d'Eilenberg-Zilber, c'est \`a dire des
transformations naturelles entre bifoncteurs. Pour $j=1,2$, soient $X_j$ un $R$-RAV
alg\'ebrique, $U_j$ un ouvert de $X_j$ et $\ph_j:S(U_j,\U_j,R)\to S(X_j,R)$ un morphisme
admissible. Le produit $X_1\ti X_2$ est encore un $R$-RAV alg\'ebrique ([4], th\'eor\`eme
4) et $\ph_0=\Psi\circ(\ph_1\otimes\ph_2)\circ\Phi$ est un morphisme admissible de
$S(U_1\ti U_2,\U_1\ti\U_2,R)$ dans $S(X_1\ti X_2,R)$. Avec ces notations, nous avons le
r\'esultat suivant:

\begin{theo}$\ind(\ph_0,U_1\ti U_2)=\ind(\ph_1,U_1).\ind(\ph_2,U_2)$. \end{theo}

La d\'emonstration de ce th\'eor\`eme n\'ecessite un chemin d\'etourn\'e, l'utilisation
directe de la d\'efinition de l'indice \'etant ici mal commode:

1) Soit $\ph:S(U,\U,R)\to S(X,R)$ un morphisme admissible, o\`u $X$ est un $R$-RAV
alg\'ebrique. Si $\ph$ est de type fini, alors, avec les notations utilis\'ees plus haut,
la restriction $\varkappa$ de $Sd_2\circ\ph\circ Sd_1\circ\ta_A$ \`a $S(X,\W,R)$ est un
endomorphisme de type fini du module gradu\'e $S(X,\W,R)$ tel que $\La(\varkappa)=
\ind(\ph,U)$ ([5], lemme 3). Cela permet de d\'emontrer le th\'eor\`eme 5 dans le cas
particulier o\`u $\ph_1$ et $\ph_2$ sont de type fini, en utilisant le fait que le nombre
de Lefschetz d'un produit tensoriel de morphismes de type fini est le produit de leurs
nombres de Lefschetz.

2) Plongeons $X_i$ comme ferm\'e dans un espace norm\'e $E_i$. D'apr\`es le th\'eor\`eme
1, $X_i$ est un r\'etracte de voisinage alg\'ebrique de $E_i$, ce qui entra\^\i ne
l'existence d'un ouvert $\wh U_i$ de $E_i$ tel que $\wh U_i\cap X_i=U_i$, d'un
recouvrement ouvert $\wh\U_i$ de $\wh U_i$ tel que $\wh\U_i|U_i$ soit plus fin que $\U_i$
et d'un morphisme de cha\^\i nes $\hat\ph_i:S(\wh U_i,\wh\U_i,R)\to S(X_i,R)$ tel que
$\hat\ph_i|S(U_i,\wh\U_i|U_i,R)=\ph_i|S(U_i,\wh\U_i|U_i,R)$. Le morphisme $\hat\ph_i$ est
admissible et il r\'esulte du lemme 1(ii) et de la propri\'et\'e de contraction que
$\ind(\ph_i,U_i)=\ind(\hat\ph_i,\wh U_i)$, ce qui permet de r\'eduire la d\'emonstration
du th\'eor\`eme 5 au cas o\`u les $X_i$ sont des espaces norm\'es.

3) Soient $U$ un ouvert d'un espace norm\'e $X$ et $\ph:S(U,\U,R)\to S(X,R)$ un morphisme
admissible. Si $U'$ est un ouvert de $X$ v\'erifiant $\Fix\ph\su U'$ et $\ov U'\su U$, il
existe un morphisme admissible $\ph':S(U',\U|U',R)\to S(X,R)$ tel que $\im\ph$ soit de
type fini et une homotopie admissible $\chi:S(U',\U|U',R)\to S(X,R)$ entre $\ph'$ et la
restriction de $\ph$, ce qui ach\`eve de ramener la d\'emonstra\-tion au cas particulier
1).

\medskip
Soit $U$ un ouvert d'un espace topologique $X$. Une fonction continue $f:U\to X$ est dite
admissible si elle est compacte et si $\Fix f$ est un sous-ensemble compact de $U$. De
m\^eme, une homotopie $F:U\ti I\to X$ est dite admissible si elle est compacte et si
$\bigcup_{t\in I}\Fix F_t$ est un sous-ensemble compact de $U$. Le morphisme $f_\#:S(U,R)
\to S(X,R)$ induit par une fonction admissible est admissible donc, si $X$ est un $R$-RAV
alg\'ebrique, nous pouvons d\'efinir l'indice de point fixe de $f$ en posant
$\ind(f,U)=\ind(f_\#,U)$. Cet indice a toutes les propri\'et\'es habituelles d'un indice
de point fixe.

\begin{theo}L'indice de point fixe ainsi d\'efini a les propri\'et\'es suivantes:
\begin{itemize}
\item[(I)] {\em (Normalisation)} Si $U=X$, alors $\ind(f,U)=\La(f)$.

\item[(II)] {\em (Additivit\'e)} Si $V_1$, $V_2$ sont deux ouverts disjoints contenus
dans $U$ tels que $\Fix f\su V_1\cup V_2$, alors $\ind(f,U)=\ind(f,V_1)+\ind(f,V_2)$.

\item[(III)] {\em (Existence)} Si $\ind(f,U)\ne0$, alors $f$ a un point fixe.

\item[(IV)] {\em (Homotopie)} Si $F:U\ti I\to X$ est une homotopie admissible, alors
$\ind(F_0,U)=\ind(F_1,U)$.

\item[(V)] {\em (Contraction)} Soit $Y$ un sous-espace de $X$ qui est aussi un $R$-RAV
alg\'ebrique. Si $f(U)$ est contenu dans un compact de $Y$ et si $f_Y:U\cap Y\to Y$ est
induite par $f$, alors $\ind(f,U)=\ind(f_Y,U\cap Y)$.

\item[(VI)] {\em (Multiplicativit\'e)} Si $f_1:U_1\to X_1$ et $f_2:U_2\to X_2$ sont
admissibles, o\`u $U_1$ et $U_2$ sont des ouverts des $R$-RAV alg\'ebriques $X_1$ et
$X_2$, alors $\ind(f_1\ti f_2,U_1\ti U_2)=\ind(f_1,U_1).\ind(f_2,U_2)$.

\item[(VII)] {\em (Commutativit\'e)} Soient $U_0$ et $U_1$ des ouverts des $R$-RAV
alg\'ebri\-ques $X_0$ et $X_1$, et soient $f_0:U_0\to X_1$ et $f_1:U_1\to X_0$ des
fonctions continues. Si les fonctions $f_1f_0:U'_0=f_0^{-1}(U_1)\to X_0$ et $f_0f_1:U'_1=
f_1^{-1}(U_0)\to X_1$ sont admissibles et si $f_1$ est compacte, alors $\ind(f_1f_0,U'_0)
=\ind(f_0f_1,U'_1)$.
\end{itemize}\end{theo}

Les propri\'et\'es (I)-(VI) r\'esultent imm\'ediatement des th\'eor\`emes 4 et 5. Par
contre, la commutativit\'e, qui n'a aucun analogue simple et g\'en\'eral pour les
morphismes de cha\^\i nes, a une d\'emonstration plus d\'elicate. Notons que,
contrairement aux constructions habituelles de l'indice, la commutativit\'e ne joue aucun
r\^ole dans notre approche de l'indice; seul est utilis\'e le cas particulier (V), ou
plut\^ot son analogue pour les morphismes de cha\^\i nes.

Une autre propri\'et\'e importante de l'indice, le th\'eor\`eme $\mod p$, s'\'etend aux
RAV alg\'ebriques. Si $U$ est un ouvert de $X$, $f:U\to X$ une fonction continue et
$m\ge1$ un entier, nous notons $f^m$ la compos\'ee de $m$ copies de $f$; cette fonction
est d\'efinie sur un sous-ensemble de $U$.

\begin{theo}Soient $X$ un $\Bbb Z$-RAV alg\'ebrique, $U$ un ouvert de $X$, $f:U\to X$ une
fonction compacte et $U'$ un ouvert contenu dans $U$ sur lequel $f^m$ est d\'efinie, o\`u
$m=p^k$ avec $p$ premier. Supposons que $S=\{x\in U'\,|\,f^m(x)=x\}$ soit compact et que
$f(S)\su S$. Alors $\ind(f^m,U')\equiv\ind(f,U')\mod p$.
\end{theo}

D'apr\`es le th\'eor\`eme 3, l'indice d\'ecrit dans cette section s'applique en
particulier \`a tous les sous-ensembles convexes m\'etrisables des e.v.t. et \`a tous les
groupes m\'etrisables localement contractiles. A titre d'illustration, le r\'esultat
suivant, qui g\'en\'eralise un th\'eor\`eme bien connu de Borsuk, est prouv\'e dans [5].

\begin{theo}Soit $G$ un groupe topologique m\'etrisable contractile tel que la fonction
$x\mapsto x^2$ soit un hom\'eomorphisme de $G$ sur $G$. Soient $U$ un voisinage ouvert
sym\'etrique de l'\'el\'ement neutre de $G$ et $f:\ov U\to G$ une fonction continue
compacte sans point fixe sur $\ov U\sm U$. Si $f(x^{-1})=f(x)^{-1}$ pour tout $x\in\ov
U\sm U$, alors $\ind(f,U)$ est impair.
\end{theo}

\section{Fonctions multivoques et morphismes de cha\^\i nes}

Toutes les fonctions multivoques consid\'er\'ees dans cet article seront suppos\'ees
faire correspondre \`a tout point $x$ d'un espace topologique $X$ un sous-ensemble
ferm\'e non vide $F(x)$ d'un espace topologique $Y$. La fonction $F$ est dite
semi-continue sup\'erieurement, ou s.c.s., si, pour tout ouvert $U$ de $Y$, l'ensemble
des $x\in X$ tels que $F(x)\su U$ est ouvert dans $X$, et $F$ est dite compacte si son
image est contenue dans un compact de $Y$.

Tous les espaces m\'etrisables consid\'er\'es dans cette section seront suppos\'es munis
d'une distance, arbitraire mais fix\'ee et not\'ee $d$, et si $A$ est un sous-ensemble
d'un tel espace $X$, nous poserons $B(A,\ep)=\{x\in X\,|\,d(x,A)<\ep\}$.

Nous n'utilserons dans cette section que des coefficients rationnels, bien que la
d\'efinition ci-apr\`es se g\'en\'eralise \`a des coefficients arbitraires.

\begin{defi}Soient $X$, $Y$ des espaces m\'etrisables et $F:X\mt Y$ une fonction
multivoque compacte s.c.s.. Nous dirons que $F$ est approximable par des morphismes de
cha\^\i nes compacts si, pour tout recouvrement ouvert $\U$ de $X$ et tout $\ep>0$, il
existe un recouvrement ouvert $\V$ de $X$ et un morphisme de cha\^\i nes compact
$\ph:S(X,\V,\Bbb Q)\to S(Y,\Bbb Q)$ tels que, pour tout $V\in\V$, il existe $U\in\U$
contenant $V$ et v\'erifiant $\ph(S(V,\Bbb Q))\su S(B(F(U),\ep),\Bbb Q)$.
\end{defi}

La compacit\'e de $F$ garantit que cette d\'efinition ne d\'epend pas de la distance
utilis\'ee sur $Y$.

Soient $X$ un $\Bbb Q$-RAV alg\'ebrique et $F:X\mt X$ une fonction multivoque compacte
s.c.s. approximable par des morphismes de cha\^\i nes compacts. Il est possible
d'associer \`a $F$ un \og nombre de Lefschetz\fg $\La(F)$, mais cette d\'efinition est un
peu inhabituelle. En effet, soient $\U$ un recouvrement ouvert de $X$ et $\ep>0$, et
supposons que, pour $i=1,2$, le recouvrement ouvert $\V_i$ et le morphisme de cha\^\i nes
compact $\ph_i:S(X,\V_i,\Bbb Q)\to S(X,\Bbb Q)$ v\'erifient les conditions de la
d\'efinition 4 relativement \`a $\U$ et $\ep$. Soit $\V$ le recouvrement ouvert de $X$
form\'e des ensembles $V_1\cap V_2$ o\`u $V_1$ appartient \`a $\V_1$ et $V_2$ \`a $\V_2$,
et soit $\ph'_i$ la restriction de $\ph_i$ \`a $S(X,\V,\Bbb Q)$, de sorte que
$\La(\ph'_{i*})=\La(\ph_{i*})$. Pour tout rationnel $q$, soit
$\ph_q=q\ph'_1+(1-q)\ph'_2$; c'est un morphisme de cha\^\i nes compact de $S(X,\V,\Bbb
Q)$ dans $S(X,\Bbb Q)$ (qui, comme nous l'imposons, pr\'eserve l'augmentation). Pour
$V=V_1\cap V_2\in \V$, il existe, pour $i=1,2$, $U_i\in\U$ contenant $V_i$ et tel que
$\ph_i(S(V_i,\Bbb Q))\su S(B(F(U_i),\ep),\Bbb Q)$. Alors $\ph_q(S(V,\Bbb Q))$ est contenu
dans
$$
S(B(F(U_1),\ep),\Bbb Q)+S(B(F(U_2),\ep),\Bbb Q)\su S(B(F(U_1\cup U_2),\ep),\Bbb Q),
$$
ce qui montre que $\V$ et $\ph_q$ v\'erifient les conditions de la d\'efinition 4
relativement \`a $\St(\U)$ et $\ep$. Si $\La(\ph_{1*})\ne\La(\ph_{2*})$, alors
$$
\La(\ph_{q*})=q\La(\ph_{1*})+(1-q)\La(\ph_{2*})
$$
est un rationnel arbitraire. Deux cas sont donc possibles:
\begin{enumerate}
\item[(1)]Il existe un recouvrement ouvert $\U_0$, un $\ep_0>0$ et un rationnel $q_0$
tels que, pour tout recouvrement ouvert $\U$ plus fin que $\U_0$ et tout $0<\ep<\ep_0$,
si le recouvrement ouvert $\V$ et le morphisme de cha\^\i nes compact $\ph:S(X,\V,\Bbb
Q)\to S(X,\Bbb Q)$ v\'erifient les conditions de la d\'efinition 4 relativement \`a $\U$
et $\ep$, alors $\La(\ph_*)=q_0$. Nous posons alors $\La(F)=q_0$.

\item[(2)]Pour tout recouvrement ouvert $\U$, tout $\ep>0$ et tout rationnel $q$, il
existe un recouvrement ouvert $\V$ et un morphisme de cha\^\i nes compact
$\ph:S(X,\V,\Bbb Q)\to S(X,\Bbb Q)$ v\'erifiant les conditions de la d\'efinition 4
relativement \`a $\U$ et $\ep$ et tels que $\La(\ph_*)=q$. Nous posons alors $\La(F)=\Bbb
Q$. \end{enumerate}

Un exemple du deuxi\`eme cas est fourni par la fonction $F:S^2\mt S^2$ d\'efinie par
$F(x)=\{x,-x\}$. Les fonctions $f(x)=x$ et $g(x)=-x$ sont des s\'elections continues de
$F$, donc, pour tout recouvrement ouvert $\U$ de $S^2$, $\V=\U$ et la restriction $\ph$
de $f_*$ ou $g_*$ \`a $S(S^2,\U,\Bbb Q)$ v\'erifient les conditions de la d\'efinition 4
relativement \`a $\U$ et \`a tout $\ep>0$, mais $\La(f_*)=2\ne0=\La(g_*)$.

\begin{theo}Soient $X$ un $\Bbb Q$-RAV alg\'ebrique et $F:X\mt X$ une fonction multivoque
compacte s.c.s. approximable par des morphismes de cha\^\i nes compacts. Si $\La(F)\ne0$,
alors $F$ a un point fixe. En particulier, $F$ a un point fixe si $\wt H(X,\Bbb Q)=0$.
\end{theo}

\begin{proof}Si $F$ n'a pas de point fixe, il existe $\ep>0$ tel que $d(x,F(x))>3\ep$ pour
tout $x\in X$. Soit $\U$ le recouvrement de $X$ par les boules ouvertes $B(x,\ep)$, $x\in
X$. Soient $\V$ et $\ph:S(X,\V,\Bbb Q)\to S(X,\Bbb Q)$ un recouvrement ouvert et un
morphisme de cha\^\i nes compact v\'erifiant les conditions de la d\'efinition 4
relativement \`a $\U$ et $\ep$. Pour tout $V\in\V$, il existe $x\in X$ tel que $B(x,\ep)$
contienne $V$ et que $S(B(F(B(x,\ep)),\ep),\Bbb Q)$contienne $\ph(S(V,\Bbb Q))$. Pour
toute cha\^\i ne $c\in S(V,\Bbb Q)$ et tout point $a$ de $||\ph(c)||$, il existe $y\in
B(x,\ep)$ et $z\in F(y)$ tels que $d(a,z)<\ep$. Nous avons alors
$$
d(a,x)\ge d(y,z)-d(a,z)-d(x,y)>3\ep-\ep-\ep=\ep,
$$
ce qui montre que $||\ph(c)||\cap V\su ||\ph(c)||\cap B(x,\ep)=\vi$. L'ouvert $V$ ne
contient donc aucun point fixe de $\ph$ et, ceci \'etant vrai pour tout $V\in\V$, $\ph$
n'a pas de point fixe, donc $\La(\ph)=0$ d'apr\`es le th\'eor\`eme 4. Nous sommes donc
dans le cas (1) et $\La(F)=0$.\end{proof}

Les deux th\'eor\`emes suivants donnent des exemples de classes auxquelles s'applique le
r\'esultat pr\'ec\'edent. Un compact est dit $\Bbb Q$-acyclique si son homologie de \v
Cech r\'eduite \`a coefficients rationnels est triviale.

\begin{theo}Soient $X$, $Y$ des espaces m\'etrisables et $F:X\mt Y$ une fonction
multivoque compacte s.c.s. telle que $F(x)$ soit $\Bbb Q$-acyclique pour tout $x\in X$.
Si $Y$ est un $\Bbb Q$-RAV alg\'ebrique, alors $F$ est approximable par des morphismes de
cha\^\i nes compacts. \end{theo}

Une fonction multivoque $F:X\mt Y$ est dite continue si elle est continue lorsque
l'ensemble des ferm\'es de $Y$ est muni de la topologie de Vietoris.

\begin{theo}Soient $X$, $Y$ des espaces m\'etrisables, $q>1$ un entier et $F:X\mt Y$ une
fonction multivoque compacte continue telle que, pour tout $x\in X$, $F(x)$ ait une ou
$q$ composantes, chacune de ces composantes \'etant $\Bbb Q$-acyclique. Si $Y$ est un
$\Bbb Q$-RAV alg\'ebrique, alors $F$ est approximable par des morphismes de cha\^\i nes
compacts. \end{theo}

La preuve du th\'eor\`eme 10 est implicite dans la d\'emonstration du th\'eor\`eme 7 de
[4], et une modification de cet argument prouve le th\'eor\`eme 11 (cette modification
est un peu inhabituelle; en particulier, avec les notations de [4], si $\si$ est un
$0$-simplexe de $X$ d'image $x$ et $n\ge1$ est un entier, $\mu(\si_n)$ est soit un
$0$-simplexe de $X$, soit une combinaison lin\'eaire $\frac1q(x_\si^1+\dots+x_\si^q)$ de
$0$-simplexes distincts selon que $F(x)$ a une ou $q$ composantes).

\section{Construction universelle et points fixes}

Soit $X$ un compact. Nous construirons inductivement une suite croissante de compacts
$\{U_n(X)\}$ et des fonctions continues $\la_n:U_{n-1}(X)\ti U_{n-1}(X)\ti I\to U_n(X)$
v\'erifiant $\la_n(x,y,0)=x$, $\la_n(x,y,1)=y$ et $\la_n(x,x,t)=x$ quels que soient $x$,
$y$ dans $U_{n-1}(X)$ et $t$ dans $I$. Posons $U_{-1}(X)=\vi$, $U_0(X)=X$, et soit
$\la_0$ la fonction vide. Supposant $U_n(X)$ et $\la_n$ construits, soit $\R$ la relation
d'\'equivalence sur $U_n(X)\ti U_n(X)\ti I$ qui, pour tout $z\in U_n(X)$, identifie
$(\{z\}\ti U_n(X)\ti\{0\})\cup(U_n(X)\ti\{z\}\ti\{1\})\cup(\{z\}\ti\{z\}\ti I)$\`a un
point et qui identifie aussi $(x,y,t)$ \`a $(\la_n(x,y,t),\la_n(x,y,t),0)$ quels que
soient $x$, $y$ dans $U_{n-1}(X)$ et $t$ dans $I$. La relation $\R$ est ferm\'ee, donc le
quotient de $U_n(X)\ti U_n(X)\ti I$ par $\R$ est un compact, que nous appelons
$U_{n+1}(X)$. Nous notons $[x,y,t]$ l'image du point $(x,y,t)$ dans $U_{n+1}(X)$. La
fonction $z\mapsto [z,z,0]$ est un plongement de $U_n(X)$ dans $U_{n+1}(X)$ par lequel
nous identifions $U_n(X)$ \`a un sous-espace de $U_{n+1}(X)$, et la fonction $\la_{n+1}:
U_n(X)\ti U_n(X)\ti I\to U_{n+1}(X)$ d\'efinie par $\la_{n+1}(x,y,t)=[x,y,t]$ est
continue, a les propri\'et\'es souhait\'ees, et la d\'efinition de $\R$ garantit que
$\la_{n+1}| U_n(X)\ti U_n(X)\ti I=\la_n$. Nous munissons $U(X)=\bigcup_{n=1}^\ii U_n(X)$
de la topologie limite inductive de la suite croissante de compacts $\{U_n(X)\}$, et nous
pouvons d\'efinir une fonction $\la_X:U(X)\ti U(X)\ti I\to U(X)$ par $\la_X|U_n(X)\ti
U_n(X)\ti I=\la_n$ pour tout $n$. La fonction $\la_X$ est continue et fait de $U(X)$ un
espace uniform\'ement contractile. Cette construction a la propri\'et\'e importante
suivante ([6], lemme 1).

\begin{lem}Soient $X$ un compact, $Y$ un espace topologique et $u:X\to Y$ une fonction
continue.

(i) Si $Y$ est UC, il existe une fonction continue $v:U(X)\to Y$ prolongeant $u$.

(ii) Si $Y$ est ULC, il existe un voisinage $O$ de $X$ dans $U(X)$ et une fonction
continue $w:O\to Y$ prolongeant $u$. \end{lem}

Soient $Y$ un espace ULC et $X$ un compact de $Y$. Le lemme 2 nous fournit un voisinage
ouvert $O$ de $X$ dans $U(X)$ et une fonction continue $r:O\to Y$ telle que $r(x)=x$ pour
tout $x\in X$. Si $f:Y\to Y$ est une fonction continue dont l'image est contenue dans
$X$, nous d\'efinissons une fonction continue $g:O\to O$ en posant $g(x)=f(r(x))$ pour
tout $x\in O$. Les fonctions $f$ et $g$ ont les m\^emes points fixes et si $f_*$  et
$g_*$ sont de type fini (ce qui est en fait toujours le cas) et si l'anneau de
coefficients est principal, alors $\La(f)=\La(g)$ (car si $f':Y\to O$ est la fonction
induite par $f$, alors $f=r\circ f'$ et $g=f'\circ r$). Si $F:Y\mt Y$ est une fonction
multivoque dont l'image est contenue dans $X$, nous d\'efinissons une fonction multivoque
$G:O\mt O$ par $G(x)=F(r(x))$. Les fonctions $F$ et $G$ ont les m\^emes points fixes et
$G$ est s.c.s. (resp. continue) si $F$ est s.c.s. (resp. continue). En outre, toute
hypoth\`ese v\'erifi\'ee par les ensembles $F(y)$, $y\in Y$, est automatiquement
v\'erifi\'ee par les ensembles $G(x)$, $x\in O$. Enfin, soient $U$ un ouvert de $Y$ et
$f:U\to Y$ une fonction dont l'image est contenue dans $X$. Si $V=r^{-1}(U)$ et si
$g:V\to O$ est d\'efinie par $g(x)=f(r(x))$, alors $g$ est admissible si $f$ l'est. Nous
pouvons ainsi r\'eduire tout probl\`eme de point fixe pour les espaces ULC au cas
particulier des fonctions d\'efinies sur des ouverts de $U(X)$ et \`a valeurs dans $X$.

Pour tout compact $X$, soit $\T(X)$ l'ensemble des topologies m\'etrisables $\tau$ sur
$U(X)$ qui sont moins fines que la topologie limite inductive et telles que $\la_X:
(U(X),\tau)\ti(U(X),\tau)\ti I\to(U(X),\tau)$ soit continue. Chaque \'el\'ement de
$\T(X)$ induit sur $X$ sa topologie d'origine, donc $\T(X)=\vi$ si $X$ n'est pas
m\'etrisable. Lorsque le compact $X$ est m\'etrisable, les \'el\'ements de $\T(X)$
permettent de ramener les probl\`emes de point fixe pour les fonctions d'un ouvert de
$U(X)$ \`a valeurs dans $X$ \`a des questions analogues pour les espaces ULC
m\'etrisables gr\^ace au lemme suivant.

\begin{lem}Soient $X$ et $C$ des compacts m\'etrisables et $V$ un ouvert de $U(X)$.

(i) Si $f:V\to C$ est une fonction continue, il existe $\tau\in\T(X)$ telle que
$f:(V,\tau)\to C$ soit continue.

(ii) Si $F:V\mt C$ est une fonction multivoque s.c.s., il existe $\tau\in\T(X)$ telle que
$F:(V,\tau)\mt C$ soit s.c.s.. \end{lem}

L'affirmation (i) est contenue dans le lemme 2 de [6]. L'affirmation (ii) est le lemme 1
de [3] quand $C=X$ et $V=U(X)$, mais la d\'emonstration de ce lemme s'applique au cas
g\'en\'eral.

Notons aussi que le passage de la topologie libre \`a une topologie m\'etrisable
$\tau\in\T(X)$ pour laquelle une fonction $f$ d\'efinie sur un ouvert de $U(X)$ reste
continue ne change pas les propri\'et\'es homologiques de $f$. Cela r\'esulte du fait
suivant, qui est un cas particulier du lemme 4 de [6].

\begin{lem}Soient $X$ un compact m\'etrisable et $\tau\in\T(X)$. Pout tout sous-ensemble
$\tau$-ouvert $V$ de $U(X)$, l'identit\'e de $(V,\tau)$ dans $V$ est une \'equivalence
homotopique faible.  \end{lem}

\section{Utilisation des $\om$-syst\`emes projectifs}

Nous expliquerons dans cette section comment l'utilisation de syst\`emes projectifs
particuliers permet de passer des compacts m\'etrisables aux compacts arbitraires dans
l'\'etude des points fixes de fonctions de $U(X)$ dans $X$.

Un syst\`eme projectif d'espaces topologiques sur un ensemble ordonn\'e filtrant $A$ sera
not\'e $\BS=(X_\al,p_\al^\bt,A)$ ( les $X_\al$ sont des espaces topologiques et, pour
$\al\le\bt$, $p_\al^\bt:X_\bt\to X_\al$ est continue ). Nous notons $\lim\BS$ la limite
projective de ce syst\`eme et $p_\al$ la projection de $\lim\BS$ dans $X_\al$. Si $B$ est
un sous-ensemble filtrant de $A$, nous notons $\BS|B=(X_\al,p_\al^\bt,B)$ le syst\`eme
obtenu en restreignant l'ensemble des indices \`a $B$. Si $B$ est cofinal dans $A$, nous
identifions naturellement $\lim(\BS|B)$ \`a $\lim\BS$. Si $\BS_1=(X_\al,p_\al^\bt,A)$ et
$\BS_2=(Y_\al,q_\al^\bt,A)$ sont deux syst\`emes projectifs ayant le m\^eme ensemble
d'indices, un morphisme de $\BS_1$ dans $\BS_2$ est une famille de fonctions continues
$f_\al:X_\al\to Y_\al$, $\al\in A$, telle que $f_\al p_\al^\bt=q_\al^\bt f_\bt$ pour
$\al\le\bt$. Un tel morphisme induit une fonction continue $\lim(f_\al)$ de $\lim\BS_1$
dans $\lim\BS_2$.

Un syst\`eme projectif $\BS=(X_\al,p_\al^\bt,A)$ est appel\'e un $\om$-syst\`eme
projectif s'il v\'erifie les conditions suivantes:

\begin{itemize}
\item[($*$)]Tout sous-ensemble d\'enombrable totalement ordonn\'e de $A$ a une borne
sup\'erieure.

\item[($**$)]Pour tout sous-ensemble totalement ordonn\'e $B$ de $A$ admettant une borne
sup\'erieure $\bt$, la fonction $\lim_{\al\in B}p_\al^\bt:X_\bt\to\lim(\BS|B)$ est un
hom\'eomorphisme.

\item[($***$)]Chaque $X_\al$ a une base d\'enombrable.
\end{itemize}

Un sous-ensemble $B$ d'un ensemble filtrant $A$ est dit ferm\'e si, pour tout
sous-ensemble totalement ordonn\'e $C\su B$ admettant une borne sup\'erieure dans $A$,
cette borne sup\'erieure appartient \`a $B$. Si $\BS=(X_\al,p_\al^\bt,A)$ est un
$\om$-syst\`eme projectif et si $B$ est un sous-ensemble cofinal et ferm\'e de $A$, alors
$\BS|B$ est aussi un $\om$-syst\`eme projectif.

La propri\'et\'e suivante des $\om$-syst\`emes projectifs, qui est un cas particulier du
th\'eor\`eme 3.1.9 de [7], joue un r\^ole essentiel dans nos arguments.

\begin{lem}Soient $\BS_1=(X_\al,p_\al^\bt,A)$ et $\BS_2=(Y_\al,q_\al^\bt,A)$ deux
$\om$-syst\`emes projectifs de compacts sur un m\^eme ensemble d'indices $A$, et soit $f$
une fonction continue de $\lim\BS_1$ dans $\lim\BS_2$. Si, pour tout $\al\in A$, la
projection $p_\al:\lim\BS_1\to X_\al$ est surjective, alors il existe un sous-ensemble
cofinal ferm\'e $B$ de $A$ et un morphisme de $\BS_1|B$ dans $\BS_2|B$ dont la limite est
$f$.  \end{lem}

Tout compact non m\'etrisable $X$ peut \^etre repr\'esent\'e comme limite d'un
$\om$-syst\`eme projectif $\BS=(X_\al,p_\al^\bt)$ tel que les projections $p_\al:X\to
X_\al$ soient surjectives. Un tel syst\`eme peut se construire comme suit. Plongeons $X$
dans un cube de Tychonoff $I^K$. Soit $A$ l'ensemble des sous-ensembles d\'enombrables de
$K$. Pour $\al$, $\bt$ dans $A$, posons $\al\le\bt$ si $\al$ est contenu dans $\bt$. Pour
$\al\in A$, soit $X_\al$ la projection de $X$ sur le sous-cube $I^\al$ et, pour
$\al\le\bt$, soit $p_\al^\bt:X_\bt\to X_\al$ la restriction \`a $X_\bt$ de la projection
de $I^\bt$ sur $I^\al$.

Dans toute la suite de cette section, $\BS=(X_\al,p_\al^\bt,A)$ d\'esignera un
$\om$-syst\`eme de compacts (m\'etrisables d'apr\`es ($***$)) de limite $X$ tel que
toutes les projections $p_\al$ soient surjectives. Les constructions $U_n(\,.\,)$ et
$U(\,.\,)$ d\'ecrites \`a la section pr\'ec\'edente sont des foncteurs. Nous noterons
$p_{\al,n}^\bt:U_n(X_\bt)\to U_n(X_\al)$, $p_{\al,n}:U_n(X)\to U_n(X_\al)$, $\td
p_\al^\bt:U(X_\bt)\to U(X_\al)$ et $\td p_\al:U(X)\to U(X_\al)$ les prolongements
naturels des fonctions $p_\al^\bt$ et $p_\al$, et nous poserons
$U_n(\BS)=(U_n(X_\al),p_{\al,n}^\bt,A)$ et $U(\BS)=(U(X_\al),\td p_\al^\bt,A)$. Les
foncteurs $U_n$, $n\ge1$, commutent aux limites projectives, mais $U$ ne commute pas
toujours aux limites projectives. Pour tout $n$, $U_n(\BS)$ est un $\om$-syst\`eme
projectif de limite $U_n(X)$, et les projections $p_{\al,n}:U_n(X)\to U_n(X_\al)$ sont
surjectives.

Soit $f:U(X)\to X$ une fonction continue. Appliquant le lemme 5 \`a la restriction de $f$
\`a $X=U_1(X)$, nous obtenons un sous-ensemble $A_1$, cofinal et ferm\'e dans $A$ et un
morphisme $(f_\al^1):U_1(\BS)|A_1\to \BS|A_1$ de limite $f|X$. Comme $U_2(\BS)|A_1$ et
$\BS|A_1$ sont des $\om$-syst\`emes projectifs, nous pouvons appliquer le lemme 5 \`a la
restriction de $f$ \`a $U_2(X)$ pour obtenir un sous-ensemble $A_2$ de $A_1$, cofinal et
ferm\'e dans $A_1$, donc aussi dans $A$, et un morphisme $(f_\al^2)_{\al\in A_2}:
U_2(\BS)|A_2\to\BS|A_2$ de limite $f|U_2(X)$. Inductivement, nous construisons ainsi une
suite d\'ecroissante $\{A_n\}$ de sous-ensembles cofinaux et ferm\'es de $A$ et des
morphismes $(f_\al^n)_{\al\in A_n}:U_n(\BS)|A_n\to \BS|A_n$ de limite $f|U_n(X)$. La
condition ($*$) garantit que $B=\bigcap_{n=1}^\ii A_n$ est cofinal dans $A$. Pour $\al\in
B$ et $n\le m$, nous avons $f_\al^n\circ p_{\al,n}=p_n\circ(f|U_n(X))=(f_\al^m|U_n(X))
\circ p_{\al,n}$; comme $p_{\al,n}$ est surjective, il en r\'esulte que $f_\al^n=
f_\al^m|U_n(X_\al)$, et nous pouvons d\'efinir une fonction continue $f_\al:U(X_\al)\to
X_\al$ par $f_\al|U_n(X_\al)=f_\al^n$ pour tout $n$. Pour $\al\le\bt$ dans $B$, nous
avons alors $f_\al\circ \td p_\al^\bt=p_\al^\bt\circ f_\bt$. Les r\'esultats des sections
pr\'ec\'edentes impliquent que l'ensemble $Z_\al$ des points fixes de $f_\al$ n'est pas
vide. Evidemment, $p_\al^\bt(Z_\bt)\su Z_\al$ pour $\al\le\bt$ dans $B$, donc
$(Z_\al,p_\al^\bt|Z_\bt,B)$ est un syst\`eme projectif de compacts non vides dont la
limite est un sous-ensemble non vide de $\lim\BS|B=X$, et chaque point de ce
sous-ensemble est un point fixe de $f$.

Plus g\'en\'eralement, cette technique a \'et\'e utilis\'ee dans [6] pour prouver le
th\'eor\`eme de Lefschetz-Hopf pour les applications compactes des espaces ULC. Nous
consid\'erons dans ce cas un voisinage ouvert $V$ de $X$ dans $U(X)$ et une fonction
continue $f:V\to V$ telle que $f(V)\su X$. Nous devons construire un sous-ensemble
cofinal $B$ de $A$ et, pour $\al\in B$, un voisinage ouvert $V_\al$ de $X_\al$ dans
$U(X_\al)$ et une fonction continue $f_\al:V_\al\to V_\al$ telle que $f_\al(V_\al)\su
X_\al$. Ces \'el\'ements doivent v\'erifier $(\td p_\al)^{-1}(V_\al)\su V$, $(\td
p_\al^\bt)^{-1}(V_\al)=V_\bt$ pour $\al\le\bt$ et $f_\al\circ(\td p_\al|(\td
p_\al)^{-1}(V_\al))=p_\al\circ(f|(\td p_\al)^{-1}(V_\al))$. Comme $X_\al$ est
m\'etrisable, $\La(f_\al)$ est d\'efini et l'ensemble $Z_\al$ des points fixes de $f_\al$
est non vide si $\La(f_\al)\ne0$. Il nous faut ensuite v\'erifier que $\La(f)$ est
d\'efini et \'egal \`a $\La(f_\al)$ pour tout $\al\in B$. Si $\La(f)\ne0$, alors
$(Z_\al,p_\al^\bt|Z_\bt,B)$ est un syst\`eme projectif de compacts non vides et tout
point de sa limite est un point fixe de $f$. La construction des $V_\al$ utilise le
r\'esultat suivant, qui est le lemme 7 de [6].

\begin{lem}Pour tout voisinage ouvert $V$ de $X$ dans $U(X)$, il existe $\al_0\in A$ et
un voisinage ouvert $V_0$ de $X_{\al_0}$ dans $U(X_{\al_0})$ tel que $(\td
p_{\al_0})^{-1}(V_0)\su V$.  \end{lem}

Soit $C(X)=X\ti I/X\ti\{1\}$ le c\^one de base $X$. Notons $[x,t]$ l'image dans $C(X)$ du
point $(x,t)$ de $X\ti I$ et $v=[x,1]$ le sommet de $C(X)$. Pour $\al\le\bt$,
d\'efinissons $\pi_\al^\bt:C(X_\bt)\to C(X_\al)$ par
$\pi_\al^\bt([x,t])=[p_\al^\bt(x),t]$. Alors $C(\BS)=(C(X_\al),\pi_\al^\bt,A)$ est un
$\om$-syst\`eme projectif de limite $C(X)$; soit $\pi_\al$ la projection de $C(X)$ sur
$C(X_\al)$. Soit $V'$ un voisinage de $X$ dans $U(X)$ dont la fermeture est contenue dans
$V$, et soit $\eta:U(X)\to I$ une fonction continue nulle sur $V'$ et \'egale \`a un sur
un voisinage de $U(X)\sm V$. La fonction $g:U(X)\to C(X)$ d\'efinie par
$$
g(y)=\begin{cases} [f(y),\eta(y)] &\text{si }y\in V \\
v &\text{sinon} \end{cases}
$$
est continue. Proc\'edant comme ci-dessus, nous pouvons construire un sous-ensemble
cofinal $B'$ de $A$ et un morphisme $(g_\al)_{\al\in B'}:U(\BS)|B'\to C(\BS)|B'$ tel que
$g_\al\circ\td p_\al=\pi_\al\circ g$. Prenons $\al_0\in A$ et un voisinage ouvert $V_0$
de $X_{\al_0}$ dans $U(X_{\al_0})$ tels que $(\td p_{\al_0})^{-1}(V_0)\su V'$. Alors $B=
\{\al\in B'\,|\,\al_0\le\al\}$ est cofinal. Pour $\al\in B$, soit $V_\al=(\td
p_{\al_0}^\al)^{-1}(V_0)$. Identifiant $X_\al$ \`a l'image de $X_\al\ti\{0\}$ dans
$C(X_\al)$, nous constatons que $g_\al(V_\al)$ est contenu dans $X_\al$ et prenons pour
$f_\al$ la fonction induite par $g_\al|V_\al$. Le fait que $\La(f)$ est d\'efini et
\'egal \`a $\La(f_\al)$ est une cons\'equence facile du r\'esultat suivant, qui est un
cas particulier du lemme 4 de [6].

\begin{lem}Pour tout $\al\in A$ et tout ouvert $O$ de $U(X_\al)$, la restriction de $\td
p_\al$ \`a $(\td p_\al)^{-1}(O)$ est une \'equivalence homotopique faible de $(\td
p_\al)^{-1}(O)$ sur $O$.  \end{lem}

La technique d\'ecrite ci-dessus s'applique aussi aux fonctions multivoques. Soit $\K(X)$
l'ensemble des compacts non vides de $X$ muni de la topologie de Vietoris. Soit
$P_\al^\bt:\K(X_\bt)\to\K(X_\al)$ le prolongement naturel de $p_\al^\bt$. Alors $\K(\BS)=
(\K(X_\al),P_\al^\bt,A)$ est un $\om$-syst\`eme projectif de limite $\K(X)$ et la
projection $P_\al:\K(X)\to\K(X_\al)$ est induite par $p_\al$.

Soit $F:U(X)\mt X$ une fonction multivoque. Si $F$ est continue, nous pouvons, comme
pr\'ec\'edemment, trouver un sous-ensemble cofinal $B$ de $A$ et, pour $\al\in B$, une
fonction continue $F_\al:U(X_\al)\to \K(X_\al)$ telle que $F_\al\circ \td
p_\al=P_\al\circ F$. Si l'ensemble $Z_\al$ des points fixes de $F_\al$ est non vide pour
tout $\al\in B$, alors $F$ a un point fixe.

Quand $F$ est seulement s.c.s., les $F_\al$ ne peuvent plus s'obtenir de la m\^eme fa\c
con, mais une autre voie est possible. D\'efinissons $F_\al:U(X_\al)\mt X_\al$ par
\begin{equation}
F_\al(x)=p_\al(F(p_{\al,n}^{-1}(x))) \text{ si $x\in U_n(X_\al)\sm U_{n-1}(X_\al)$,
$n\ge1$.}  \tag{\dag}
\end{equation}

Gr\^ace au lemme suivant ([3], lemme 4), nous pouvons continuer comme pr\'ec\'edemment.

\begin{lem}L'ensemble des $\al\in A$ tels que $F_\al$ soit s.c.s. est cofinal dans $A$.
\end{lem}

Certaines conditions de nature homologique se transmettent aux $F_\al$ pour un
sous-ensemble cofinal d'indices. Par exemple, le r\'esultat suivant, dans lequel les
$F_\al$ sont d\'efinies par la formule ($\dag$), r\'esulte du lemme 3 de [3] et de la
condition ($*$).

\begin{lem}Si $F$ est s.c.s. et si les $F(x)$ sont $\Bbb Q$-acycliques pour tout $x\in
U(X)$, alors l'ensemble des $\al\in A$ tels que $F_\al(x)$ soit $\Bbb Q$-acyclique pour
tout $x\in U(X_\al)$ est cofinal et ferm\'e dans $A$. \end{lem}

Les lemmes 9 et 8 ont \'et\'e utilis\'es dans [3] pour prouver qu'une fonction multivoque
s.c.s. compacte $F$ d'un espace uniform\'ement contractile dans lui-m\^eme a un point
fixe si tous les compacts $F(x)$ sont $\Bbb Q$-acycliques \footnote{L'article [3] a
\'et\'e \'ecrit avant la d\'ecouverte des RAV alg\'ebriques. Comme nous l'avons
mentionn\'e dans l'introduction, cette nouvelle approche a rendu inutile la correction de
l'erreur contenue dans [2]. Par cons\'equent, l'article \og Le probl\`eme de Schauder;
correction et compl\'ement\fg, mentionn\'e dans la bibliographie de [2] ne sera jamais
publi\'e. Le cas m\'etrisable de la proposition de [2] r\'esulte du lemme 1 de [2] et du
th\'eor\`eme 7 de [4], ce dernier rempla\c cant l'utilisation de l'article non
publi\'e.}.

\section{Un exemple}

La technique de d\'emonstration d\'ecrite dans les deux sections pr\'ec\'edentes est
g\'en\'erale. Pour la r\'ecapituler et l'illustrer, nous prouverons le r\'esultat
suivant.

\begin{theo}Soient $Y$ un espace UC s\'epar\'e, $q>1$ un entier et $F:Y\mt Y$ une
fonction multivoque compacte continue telle que, pour tout $x\in Y$, $F(x)$ ait une ou
$q$ composantes, chacune de ces composntes \'etant $\Bbb Q$-acyclique. Alors $F$ a un
point fixe.  \end{theo}

\begin{proof}Premi\`ere \'etape: prouver le r\'esultat quand $Y$ est m\'etrisable. Il
suffit d'appliquer les th\'eor\`emes 3, 9 et 11.

Deuxi\`eme \'etape: r\'eduction \`a un cas particulier. Soit $X$ un compact de $Y$
contenant $F(X)$, et soit $r:U(X)\to Y$ une fonction continue telle que $r(x)=x$ pour
tout $x\in X$. La fonction $G:U(X)\mt U(X)$ d\'efinie par $G(x)=F(r(x))$ v\'erifie les
m\^emes hypoth\`eses que $F$ et a les m\^emes points fixes. Le probl\`eme est donc
r\'eduit au cas particulier $Y=U(X)$ et $F(Y)\su X$.

Troisi\`eme \'etape: prouver le cas particulier lorsque le compact $X$ est m\'etrisable.
Comme $\K(X)$ est alors m\'etrisable, le lemme 3(i) nous fournit une topologie
m\'etrisable $\tau\in\T(X)$ telle que $F:(U(X),\tau)\to \K(X)$ soit continue, donc $F$ a
un point fixe d'apr\`es la premi\`ere partie de la d\'emonstration.

Quatri\`eme \'etape: passage des compacts m\'etrisables aux compacts arbitraires dans le
cas particulier. Repr\'esentons $X$ comme limite d'un $\om$-syst\`eme projectif $\BS=
(X_\al,p_\al^\bt,A)$ de compacts m\'etrisables tel que toutes les projections $p_\al:X\to
X_\al$ soient surjectives. Il nous faut trouver un sous-ensemble cofinal $B$ de $A$ et un
morphisme $(F_\al)_{\al\in B}$ de $U(\BS)|B$ dans $\K(\BS)|B$ tel que, pour tout $\al\in
B$, $F_\al\circ\td p_\al=P_\al\circ F$ et que les $F_\al$ v\'erifient les hypoth\`eses du
th\'eor\`eme. D'apr\`es la troisi\`eme \'etape, l'ensemble $Z_\al$ des points fixes de
$F_\al$ sera alors non vide, et tout point de la limite du syst\`eme projectif $(Z_\al,
p_\al^\bt|Z_\bt,B)$ sera un point fixe de $F$.

En utilisant le lemme 5 comme indiqu\'e dans la section pr\'ec\'edente, nous pouvons
trouver un sous-ensemble cofinal et ferm\'e $A_1$ de $A$ et un morphisme $(F_\al)_{\al\in
A_1}:U(\BS)|A_1\to\K(\BS)|A_1$ tel que $F_\al\circ\td p_\al=P_\al\circ F$ pour tout
$\al\in A_1$. Soit $A_2$ l'ensemble des $\al\in A_2$ tels que $F_\al(x)$ ait une ou $q$
composantes pour tout $x\in U(X_\al)$. Montrons que $A_2$ est cofinal et ferm\'e dans
$A$. Notons d'abord que si $C$ est un sous-ensemble totalement ordonn\'e de $A_1$
admettant la borne sup\'erieure $\ga$, alors la condition ($**$), appliqu\'ee \`a
$\K(\BS)|A_1$, implique que
$$
F_\ga(x)=\lim_{\bt\in C}P_\bt^\ga(F_\ga(x))=\lim_{\bt\in C}F_\bt(\td p_\bt^\ga(x)),
$$
ce qui entra\^\i ne que $A_2$ est ferm\'e. De m\^eme, nous avons
$$
F(x)=\lim_{\al\in A_1} F_\al(\td p_\al(c))
$$
pour tout $x\in U(X)$. Nous d\'eduirons la cofinalit\'e de $A_2$ du fait suivant:

\medskip\noindent
{\bf Affirmation.} {\it Pour tout $\al\in A_1$, il existe $\al^*\ge\al$ dans $A_1$ tel
que $F_{\al^*}(y)$ ait $q$ composantes pour tout $y\in U(X_{\al^*})$ tel que
$F_\al(p_\al^{\al^*}(y))$ ait plus d'une composante.}

\medskip
Pour prouver cette affirmation, il suffit de montrer que, pour tout entier $n\ge1$, il
existe $\ga_n\ge\al$ dans $A_1$ tel que $F_{\ga_n}(y)$ ait $q$ composantes pour tout
$y\in U_n(X_{\ga_n})$ tel que $F_\al(p_\al^{\ga_n}(y))$ ait plus d'une composante. En
effet, puisque $A$ est filtrant, la condition ($*$) garantit l'existence d'un $\al^*\in
A_1$ tel que $\al^*\ge\ga_n$ pour tout $n$. Si $y\in U(X_{\al^*})$ est tel que
$F_\al(p_\al^{\al^*}(y))$ ait plus d'une composante et si $n$ est tel que $y$ appartienne
\`a $U_n(X_{\al^*})$, le choix dz $\ga_n$ garantit que $F_{\ga_n}(p_{\ga_n}^{\al^*}(y))$
a $q$ composantes. Puisque
$F_{\ga_n}(p_{\ga_n}^{\al^*}(y))=P_{\ga_n}^{\al^*}(F_{\al^*}(y))$ et que $F_{\al^*}(y)$ a
au plus $q$ composantes, $F_{\al^*}(y)$ a exactement $q$ composantes.

Fixons un entier $n\ge1$, et soit $G_n$ l'ensemble des $x\in U_n(X_\al)$ tels que
$F_\al(x)$ ait plus d'une composante. Soient $x\in G_n$ et $z\in(p_{\al,n})^{-1}(x)$.
Puisque $F_\al(x)=P_\al(F(z))$ a plus d'une composante, $F(z)$ en a $q$ et, puisque
$F(z)=\lim_{\bt\in A_1}F_\bt(\td p_\bt(z))$, il existe $\bt_z\ge \al$ dans $A_1$ tel que
$F_{\bt_z}(\td p_{\bt_z}(z))$ ait $q$ composantes; soient $K_1,\dots,K_q$ ces
composantes, et soient $O_1,\dots,O_q$ des ouverts disjoints de $X_{\bt_z}$ tels que
$K_i\su O_i$ pour $1\le i\le q$. La continuit\'e de $F$ entra\^\i ne l'existence d'un
voisinage $V_z$ de $z$ dans $U_n(X)$ tel que, pour tout $z'\in V_z$, $F(z')\su
\bigcup_{i=1}^qO_i$ et $F(z')\cap O_i\ne\vi$ pour $1\le i\le q$. Le compact
$(p_{\al,n})^{-1}(x)$ peut \^etre recouvert par un nombre fini de ces ensembles $V_z$,
soient $V_{z_1},\dots,V_{z_r}$, et il existe un voisinage ouvert $W_x$ de $x$ dans
$U_n(X_\al)$ tel que $(p_{\al,n})^{-1}(W_x)\su \bigcup_{j=1}^rV_{z_j}$. Soit $\ga_x\in
A_1$ tel que $\ga_x\ge\bt_{z_j}$ pour $1\le j\le r$. Si $z$ appartient \`a
$(p_{\al,n})^{-1}(W_x)$, il existe $1\le j\le r$ tel que $z$ appartienne \`a $V_{z_j}$,
donc $F_{\bt_j}(\td p_{\bt_{z_j}}(z))=P_{\bt_{z_j}}^{\ga_x}(F_{\ga_x}(\td p_{\ga_x}(z)))$
a $q$ composantes, et l'ensemble $F_{\ga_x}(\td p_{\ga_x}(z))$ doit aussi avoir $q$
composantes.

L'ensemble $G_n$ est m\'etrisable et s\'eparable, donc peut \^etre recouvert par une
famille d\'enombrable $\{W_{x_m}\}_{m=1}^\ii$ des ouverts $W_x$, $x\in G_n$. En utilisant
($*$), nous pouvons trouver $\ga_n\in A_1$ tel que $\ga_n\ge\ga_{x_m}$ pour tout $m$.
Soit $y\in U_n(X_{\ga_n})$ tel que $F_\al(p_\al^{\ga_n}(y))$ ait plus d'une composante,
et soit $u$ tel que $W_{x_u}$ contienne $p_\al^{\ga_u}(y)$. La surjectivit\'e de
$p_{\ga_u}$ entra\^\i ne la surjectivit\'e de $p_{\ga_u,n}$; soit $z\in U_n(X)$ tel que
$p_{\ga_u,n}(z)\,(=\td p_{\ga_u}(z))=y$. Alors $F_{\ga_{x_u}}(p_{\ga_{x_u}}^{\ga_n}(y))=
F_{\ga_{x_u}}(\td p_{\ga_{x_u}}(z))$ a $q$ composantes et, comme $F_{\ga_{x_u}}
(p_{\ga_{x_u}}^{\ga_n}(z))=P_{\ga_{x_u}}^{\ga_n}(F_{\ga_n}(y))$, l'ensemble
$F_{\ga_n}(y)$ a aussi $q$ composantes.

Partant d'un \'el\'ement arbitraire $\al_0$ de $A_1$, d\'efinissons inductivement une
suite $\{\al_n\}$ d'\'el\'ements de $A_1$ par $\al_{n+1}=\al_n^*$. La borne sup\'erieure
$\bt$ de la suite $\{\al_n\}$ appartient \`a $A_1$. Si $y$ est un \'el\'ement de
$U(X_\bt)$ tel que $F_\bt(y)$ ait plus d'une composante, la relation
$F_\bt(y)=\lim_nF_{\al_n}(\td p_{\al_n}^\bt(y))$ entra\^\i ne l'existence d'un $n$ tel
que $F_{\al_n}(\td p_{\al_n}^\bt(y))$ ait plus d'une composante. Puisque
$\al_{n+1}=\al_n^*$, $F_{\al_{n+1}}(\td p_{\al_{n+1}}^\bt(y))$ a alors $q$ composantes et
l'\'egalit\'e $F_{\al_{n+1}}(\td p_{\al_{n+1}}^\bt(y))=P_{\al_{n+1}}^\bt (F_\bt(y))$
entra\^\i ne que $F_\bt(y)$ a aussi $q$ composantes. Cela montre que $\bt$ appartient \`a
$A_2$. Comme $\al_0$ est un \'el\'ement arbitraire du sous-ensemble cofinal $A_1$, la
cofinalit\'e de $A_2$ en r\'esulte.

Pour achever la d\'emonstration du th\'eor\`eme 12, il suffit de v\'erifier la
cofinalit\'e de l'ensemble $B$ des \'el\'ements $\al$ de $A_2$ tels que, pour tout $x\in
U(X_\al)$, chaque composante de $F_\al(x)$ soit $\Bbb Q$-acyclique. Soit $B_n$ l'ensemble
des $\al\in A_2$ tels que, pour tout $x\in U_n(X_\al)$, chaque composante de $F_\al(x)$
soit $\Bbb Q$-acyclique, de sorte que $B=\bigcap_{n=1}^\ii B_n$ et $B_{n+1}\su B_n$ pour
tout $n$. Comme l'intersection d'une suite d\'ecroissante de sous-ensembles cofinaux et
ferm\'es de $A$ est encore un sous-ensemble cofinal et ferm\'e, il suffit de montrer que
chaque $B_n$ est cofinal et ferm\'e dans $A$, ce qui peut se faire en adaptant les
arguments du lemme 3 de [3].
\end{proof}

\end{document}